\documentclass[a4paper,11pt,openany]{amsart}
\usepackage[centertags]{amsmath}
\usepackage{amscd}
\usepackage{amsthm}
\usepackage{amssymb}
\usepackage{array}
\usepackage{enumerate}
\usepackage{multicol}
\usepackage{centernot}
\usepackage{mathtools}
\usepackage{stmaryrd}
\usepackage{tikzsymbols}
\usepackage{textcomp}
\usepackage{parskip}
\usepackage{tikz}
\usetikzlibrary{cd}
\usepackage{tikz-cd}
\usetikzlibrary{matrix}
\usepackage{setspace}
\usepackage[utf8]{inputenc}
\usepackage[T1]{fontenc}
\usepackage{hyperref}
\usepackage{cleveref}
\usepackage[english]{babel}
\usepackage{color}
\usepackage{tikz}
\usepackage[sc]{mathpazo}
\usepackage{geometry}
\usepackage{subfiles}

\newtheorem{defn0}{Definition}[section]
\newtheorem{prop0}[defn0]{Proposition}
\newtheorem{thm0}[defn0]{Theorem}
\newtheorem{lemma0}[defn0]{Lemma}
\newtheorem{corollary0}[defn0]{Corollary}
\newtheorem{example0}[defn0]{Example}
\newtheorem{remark0}[defn0]{Remark}
\newtheorem{conjecture0}[defn0]{Conjecture}
\newtheorem{ques0}[defn0]{Question}

\newcommand{\Z}{\mathbb{Z}}
\newcommand{\Zc}{\mathcal{Z}}
\newcommand{\R}{\mathbb{R}}

\newcommand{\C}{\mathbb{C}}

\renewcommand{\ker}{\operatorname{Ker}}

\newcommand{\Hom}{\operatorname{Hom}}

\newcommand{\Homeo}{\operatorname{Homeo}}
\newcommand{\Diff}{\operatorname{Diff}}
\newcommand{\Out}{\operatorname{Out}}
\newcommand{\Aut}{\operatorname{Aut}}
\newcommand{\Inn}{\operatorname{Inn}}
\newcommand{\SO}{\operatorname{SO}}

\newcommand{\SU}{\operatorname{SU}}

\newcommand{\Gl}{\operatorname{GL}}

\newcommand{\D}{\operatorname{disc-sym}}
\newcommand{\A}{\operatorname{tor-sym}}
\newcommand{\rank}{\operatorname{rank}}
\newcommand{\Stab}{\operatorname{Stab}}
\newcommand{\PSL}{\operatorname{PSL}}

\newcommand{\Aff}{\operatorname{Aff}}
\usepackage{fancyhdr}
\begin{document}

\title{Actions of large finite groups on aspherical manifolds}
\author{Jordi Daura Serrano}
\address{Jordi Daura Serrano, Department de Màtematiques i Informàtica, Universitat de Barcelona (UB), Gran Via de les Corts Catalanes 585, 08007 Barcelona (Spain)}
\email{jordi.daura@ub.edu}

\maketitle
\pagenumbering{arabic}
\begin{abstract}
In this paper we study actions of finite groups on closed connected aspherical manifolds. Under some assumptions on the outer automorphism group of the fundamental group of a closed connected aspherical manifold $M$, we prove that the homeomorphism group of $M$ is Jordan, we bound the discrete degree of symmetry of $M$ and obtain a rigidity result, and we study the number of stabilizers of finite group actions on $M$. Thereafter, we show that closed connected aspherical locally homogeneous spaces $H\setminus G/\Gamma$ satisfy the necessary hypothesis on the outer automorphism group of the fundamental group.
\end{abstract}
\date{}

\section{Introduction}\label{sec:intro}

Given a closed topological manifold $M$, can we determine which finite groups act effectively on $M$? The answer of this question is completely out of reach in the vast majority of cases with the current tools of finite transformation groups theory. One way to modify this question to obtain a more tractable problem is to consider actions of a finite group $F$ on $M$ and to study properties not of the action of $F$ on $M$, but on the restriction to some subgroup of $F$ of index bounded by a constant $C$ only depending on $M$. Let us recall some problems that follow this philosophy from \cite{riera2023actions}, a recent survey on the topic.

We need the following definition to state the first problem:

\begin{defn0}\label{def:Jordan property}
A group $\mathcal{G}$ is said to be Jordan if there exists a constant $C$ such that every finite subgroup $F\leq \mathcal{G}$ has an abelian subgroup $A\leq F$ such that $[F:A]\leq C$. 
\end{defn0}

The name of this property is inspired by a classical theorem of Camille Jordan which states that $\Gl(n,\R)$ is Jordan. Around 30 years ago Étienne Ghys asked in a series of talks whether the diffeomorphism group of a closed connected smooth manifold $M$ is Jordan. This question has been answered affirmatively for a lot of different manifolds like closed flat manifolds, integral homology spheres, closed connected manifolds with non-zero Euler characteristics or closed connected manifolds up to dimension 3 (see \cite{IgnasiMundetiRiera2010Jtft,zimmermann2014jordan,ye2019symmetries,mundet2019finite}), while it has been shown that there are closed manifolds whose diffeomorphism group is not Jordan, like $S^2\times T^2$ (see \cite{csikos2014diffeomorphism,riera2017non}). More recently, the same question has been studied for the homeomorphism group of closed topological manifolds, extending the results obtained in the smooth case (see \cite{riera2023actions} and references therein). Note that this property does not tell us anything for finite groups $F$ acting on $M$ with $|F|\leq C$, but it becomes relevant for large enough finite groups acting effectively on $M$.

Another interesting problem (specially when $\Homeo(M)$ is Jordan) is to study the rank of finite abelian groups (defined as the smallest number of elements needed to generate the group) which can act effectively on $M$. Here again, we follow the philosophy of looking for a statement that eventually requires to replace the our group by a subgroup of bounded index. Fix a natural number $k$. Does there exist a constant $C$ such that every finite abelian group $A$ acting effectively on $M$ has an abelian subgroup $B$ such that $[A:B]\leq C$ and $\rank(B)\leq k$?. This question can be reformulated in terms of the following invariant introduced in \cite{riera2021topological}.

\begin{defn0}\label{def:discsym}
Given a manifold $M$ let 
$$\mu(M)=\{r\in\mathbb{N}:\text{ $M$ admits an effective action of $(\Z/a)^r$ for arbitrarily large $a$}\}.$$

More explicitly, $r\in \mu(M)$ if there exists an increasing sequence of natural number $\{a_i\}$ and effective group actions of $(\Z/a_i)^r$ on $M$ for each $i$. 

The discrete degree of symmetry of a manifold $M$ is
$$ \D(M)=\max(\{0\}\cup\mu(M)).$$
\end{defn0}

The above question is equivalent to asking whether $\D(M)\leq k$. By a theorem of L.N.Mann and J.C.Su (see \cite{mann1963actions}) we know that if $M$ is a closed connected manifold then $\D(M)$ is a well-defined natural number, but finding the exact value of $\D(M)$ is probably difficult in most cases. Note that  a manifold $M$ can admit group actions of abelian group of higher rank than $\D(M)$. For example, for any natural numbers $a,b$ there exists a closed surface $\Sigma_{g(a,b)}$ of genus $g(a,b)\geq 2$ such that $(\Z/a)^b$ acts freely on $\Sigma_{g(a,b)}$. On the other hand $\D(\Sigma_{g(a,b)})=0$ (see  \cite[Chapter 7]{farb2011primer}).

The definition of $\D(M)$ is analogous to the definition of the toral degree of symmetry 
$$ \A(M)=\max(\{0\}\cup\{r\in\mathbb{N}:T^r \text{ acts effectively on $M$}\})$$
studied in \cite[Chapter VII. \S2]{hsiang2012cohomology}. Recall that a classical result states that if $M$ is a closed manifold of dimension $n$ then $\A(M)\leq n$ and that $\A(M)=n$ if and only if $M\cong T^n$. However, it is not known whether $\D(M)\leq n$, and whether $\D(M)=n$ if and only if $M\cong T^n$.

The third problem is related to the number of stabilizers of a group action on a manifold.

\begin{defn0}\label{def:Stab}
Let $F$ be a finite group acting effectively on a manifold $M$. The set of stabilizer subgroups of the action of $F$ on $M$ is denoted by
$$\Stab(F,M)=\{F_x:x\in M\}.$$
\end{defn0}

It is not possible to bound $|\Stab(F,M)|$ only depending on $M$. For example, for all $n$ there is an effective action of the dihedral group $D_n$ on $S^1$ such that $|\Stab(D_n,S^1)|\geq n/2$. On the other hand, it is proven in \cite[Theorem 1.3]{csikos2021number} that for any closed connected manifold $M$ there exists a constant $C$ only depending on $M$ such that any finite $p$-group $F$ acting effectively on $M$ has a subgroup $H$ such that $[F:H]\leq C$ and $|\Stab(H,M)|\leq C$. This result was crucial to prove a generalized version of the Jordan property of homeomorphisms group of closed manifolds in \cite{csikos2022finite}, which states that there exists a constant $C$ such that every finite group $F$ acting on $M$ has a nilpotent subgroup $N$ such that $[F:N]\leq C$. It is not known if we can remove the hypothesis of $F$ being a $p$-group in \cite[Theorem 1.3]{csikos2021number}.

For the final question, we introduce the following definition:

\begin{defn0}\label{def: Minkowski property}
A group $\mathcal{G}$ is said to be Minkowski if there exists a constant $C$ such that every finite subgroup $F\leq \mathcal{G}$ fulfils that $|F|\leq C$. 
\end{defn0}  

\begin{remark0}
This name is motivated by a classical result of Hermann Minkowski which states that $\Gl(n,\Z)$ is Minkowski. The Minkowski property was studied in \cite{popov2018jordan,golota2023finite} under the name of bounded finite subgroups property.
\end{remark0}

If $M$ is a closed manifold and $\Homeo(M)$ is Minkowski then $M$ is said to be almost-asymmetric. In the particular case where $M$ does not admit any effective finite group action we say that $M$ is asymmetric. This case has been extensively studied (see \cite{puppe2007manifolds} and references therein).

In this paper we study these questions when $M$ is a closed connected aspherical manifolds, generalizing some existing results in the literature. For example, in \cite{ye2019symmetries} it is proven that $\Homeo(M)$ is Jordan when $M$ is closed flat manifold and in \cite{JoLee} a bound of the rank of elementary abelian group acting on closed solvmanifolds is given. However, it seems that these questions have not been studied in full generality on the literature. They are particularly relevant for closed connected aspherical manifolds, since tori are the only compact connected Lie groups that act effectively on them \cite[Thereom 3.1.16]{lee2010seifert}. Therefore, $\A(M)$ and $\D(M)$ become important invariants. 

We recall some notation before stating the main results of this paper. Given a group $\mathcal{G}$ we denote the automorphism group of $\mathcal{G}$ by $\Aut(\mathcal{G})$. Given $g\in\mathcal{G}$, we denote by $c_g:\mathcal{G}\longrightarrow \mathcal{G}$ the conjugation by $g$, $c_g(h)=ghg^{-1}$. The normal subgroup $\{c_g:g\in G\}\trianglelefteq \Aut(\mathcal{G})$ is denoted by $\Inn(\mathcal{G})$. Recall that $\Inn(\mathcal{G})\cong\mathcal{G}/\Zc(\mathcal{G})$, where $\Zc(\mathcal{G})$ denotes the center of $\mathcal{G}$. The outer automorphisms group of $\mathcal{G}$ is $\Out(\mathcal{G})=\Aut(\mathcal{G})/\Inn(\mathcal{G})$. 

Our first theorem is the following:

\begin{thm0}\label{main theorem1}
	Let $M$ be a closed connected $n$-dimensional aspherical manifold such that $\Zc(\pi_1(M))$ is finitely generated and $\Out(\pi_1(M))$ is Minkowski. Then:
	\begin{itemize}
		\item[1.] $\Homeo(M)$ is Jordan.
		\item[2.] $\D(M)\leq \rank\Zc(\pi_1(M))\leq n$, and $\D(M)=n$ if and only if $M$ is homeomorphic to $T^n$.
		\item[3.] If $\chi(M)\neq 0$ then $M$ is almost-asymmetric.
		\item[4.] If $\Aut(\pi_1(M))$ is Minkowski, then there exists a constant $C$ such that every finite group $F$ acting effectively on $M$ has a subgroup $H$ such that $[F:H]\leq C$ and $|\Stab(H,M)|\leq C$.
	\end{itemize}
\end{thm0}

This result is mainly a combination of several known results. Note that items $2$ and $4$ partially answer affirmatively the questions \cite[Question 3.4, Question 3.5]{riera2023actions} and \cite[Question 12.2]{riera2023actions} respectively. 

As an application of \cref{main theorem1}, we prove:

\begin{prop0}\label{aspherical Jordan cohomology}
There exists a closed connected aspherical manifold $M$ such that $\Homeo(M)$ is Jordan and $H^*(M)\cong H^*(T^2\times S^3)$.
\end{prop0}

Note that $\Homeo(T^2\times S^3)$ is not Jordan by \cite{riera2017non}. \Cref{aspherical Jordan cohomology} provides the first example of two manifolds with the same cohomology with integer coefficients such that one has Jordan homeomorphism group and the other does not.

We also obtain a rigidity result for closed connected aspherical $n$-dimensional  manifolds when $\D(M)=n-1$. Let $K$ denote the Klein bottle and $SK$ denote the only non-trivial principal $S^1$-bundle over $K$, then:

\begin{prop0}\label{submain theorem 1}
Let $M$ be a closed connected $n$-dimensional aspherical manifold such that $\Zc(\pi_1(M))$ is finitely generated and $\Out(\pi_1(M))$ is Minkowski. Assume that $\Inn\pi_1(M)$ has an element of infinite order. If $\D(M)=n-1$ then $M\cong T^{n-2}\times K$ or $M\cong T^{n-3}\times SK$
\end{prop0}

It is important to know when the hypothesis on the fundamental group of \cref{main theorem1} are fulfilled. At the moment there are no known closed aspherical manifolds where $\Zc(\pi_1(M))$ is not finitely generated. Regarding the second hypothesis, we prove: 

\begin{thm0}\label{main theorem2}
Let $\Gamma$ be a lattice in a connected Lie group $G$. Then $\Out(\Gamma)$ and $\Aut(\Gamma)$ are Minkowski.
\end{thm0}

The proof of this results uses several theorems of the theory of lattices of Lie groups, like the Mostow-Prasad-Margulis rigidity theorem, the Borel density theorem, Margulis superrigidity theorem and Margulis normal subgroup theorem, as well as the results in \cite{MalfaitWim2002Toag}, which are used to compute the outer automorphism group of a group extension.

If $G$ is a connected Lie group, $H$ is a maximal compact subgroup and $\Gamma$ is a cocompact lattice of $G$ then the closed aspherical locally homogeneous space $H\setminus G/\Gamma$ satisfies the hypothesis of \cref{main theorem1}. In particular, flat manifolds, nilmanifolds, almost-flat manifolds, solvmanifolds, infra-solvmanifolds and aspherical locally symmetric spaces satisfy the hypothesis of \cref{main theorem1}. Note that if we remove the asphericity hypothesis then closed locally homogeneous spaces do not necessarily have Jordan homeomorphism group. Indeed, $T^2\times S^2$ is homogeneous (and hence locally homogeneous) but $\Homeo(T^2\times S^2)$ is not Jordan.

\Cref{main theorem2} generalizes \cite[Theorem 1.7]{golota2023finite}, where it is proven that the outer automorphism group of a cocompact lattices on connected complex Lie groups is Minkowski. In the real case we need to be careful with the compact factors and factors isomorphic to $\PSL(2,\R)$ of the semisimple part of $G$.

Note that \cref{main theorem2} is also valid for non-cocompact lattices, although it cannot be used to deduce properties of large finite groups actions on non-compact aspherical locally homogeneous spaces, since the compactness hypothesis in \cref{main theorem1} is essential. 

To complement \cref{main theorem2}, we also prove that the bound on the discrete degree of symmetry is reached for aspherical locally homogeneous space.

\begin{thm0}\label{submain theorem2 discsym}
Let $H\setminus G/\Gamma$ be an aspherical locally homogeneous space. Then $\D(H\setminus G/\Gamma)=\rank\Zc\Gamma$.
\end{thm0}

This result is a combination of the results in \cite[Section 11.7]{lee2010seifert} on the toral degree of symmetry and the validity of the Borel conjecture for lattices of connected Lie groups \cite{kammeyer2016farrell}.  

Finally, using similar arguments used to prove \cref{main theorem2} we can also prove:

\begin{prop0}\label{hyperbolic product asymmetric}
Let $M=M_1\times \cdots\times M_m$, where $M_i$ are a closed aspherical manifolds such that $\pi_1(M_i)$ is hyperbolic and $\dim(M_i)\geq 3$. Then $\Out(\pi_1(M))$ is finite and $\Aut(\pi_1(M))$ is Minkowski.
\end{prop0}

\Cref{hyperbolic product asymmetric} together with \cref{main theorem1} and the fact that $\pi_1(M)$ is centreless implies that $M$ is almost asymmetric. 

The paper is divided as follows. In the second section we review the construction in \cite{MalfaitWim2002Toag} used to compute the outer automorphism group of a group extension. We also make some observations on the Minkowski property that will be used thereafter. In the third section we prove \cref{main theorem1} and we make some remarks on the relation between the discrete degree of symmetry and coverings of closed connected aspherical manifolds. In section $4$ we deduce \cref{main theorem2} for solvable Lie groups and give some examples and in section 5 we prove \cref{main theorem2} for semisimple Lie groups and \cref{hyperbolic product asymmetric}. Section 6 is devoted to combining the proofs of section 4 and section 5 in order to complete the proof of \cref{main theorem2} and we also prove \cref{submain theorem2 discsym}. Finally, in section 7 we prove \cref{aspherical Jordan cohomology} and in section 8 we prove \cref{submain theorem 1}.

\textbf{Acknowledgements:} First of all, I would like to thank Ignasi Mundet i Riera for all the guidance given during this project and for the extensive revision of the first draft of this paper. I would also like to thank László Pyber and Endre Szabó for fruitful discussions, which helped to improve many results of this paper, and for their hospitality during my stay at Álfred Rényi Institute of Mathematics. This work was partially supported by the grant PID2019-104047GB-I00 from the Spanish Ministry of Science and Innovation and the Departament de Recerca i Universitats de la Generalitat de Catalunya (2021 SGR 00697).


\section{Preliminaries on outer automorphism group}\label{sec:preliminares}

The aim of this section is to briefly explain the constructions in \cite{MalfaitWim2002Toag} used to compute the outer automorphism group of a group extension. 

Let
\[\begin{tikzcd}
1\ar{r}{}&K\ar{r}{}&G\ar{r}{p}&Q\ar{r}{}&1
\end{tikzcd}\]
be a short exact sequence of groups. The extension is determined by the morphism $\psi:Q\longrightarrow\Out(K)$ (called the abstract kernel) such that $\psi(q)$ is the class of the conjugation $c_{\sigma(q)|K}:K\longrightarrow K$, where $\sigma:Q\longrightarrow G$ is a set-theoretic section, and a 2-cocycle $c\in H^2_\psi(Q,\Zc K)$. Let $\Aut(G,K)=\{f\in\Aut(G):f(K)=K\}$ and $\Out(G,K)=\Aut(G,K)/\Inn G$. Finally, recall that if $H$ is a subgroup of $G$, the centralizer $C_G(H)$ is $\{g\in G:c_h(g)=g\text{ for all }h\in H\}$.

\begin{thm0}\label{out ses}\cite{MalfaitWim2002Toag}
There exist short exact sequences
\[\begin{tikzcd}
	1\ar{r}{}&\mathcal{K}\ar{r}{}&\Out(G,K)\ar{r}{}&L_1\ar{r}{}&1,
\end{tikzcd}\]
\[\begin{tikzcd}
	1\ar{r}{}&\overline{H}_\psi^1(Q,\Zc K)\ar{r}{}&\mathcal{K}\ar{r}{}&L_2\ar{r}{}&1
\end{tikzcd}\]
and \[\begin{tikzcd}
	1\ar{r}{}&\overline{B}_\psi^1(Q,\Zc K)\ar{r}{}&Z_\psi^1(Q,\Zc K)\ar{r}{}&\overline{H}_\psi^1(Q,\Zc K)\ar{r}{}&1
\end{tikzcd}\]
such that $L_1\leq \Out(Q)$, $L_2\leq C_{\Out(K)}\psi(Q)/\psi(\Zc Q)$ and $B_\psi^1(Q,\Zc K)\leq\overline{B}_\psi^1(Q,\Zc K)$, thus we have a surjective map $H^1_\psi(Q,\Zc K)\longrightarrow \overline{H}_\psi^1(Q,\Zc K)$. 

We have an isomorphism $\overline{B}_\psi^1(Q,\Zc K)\cong (p^{-1}(\Zc(Q))\cap C_GK)/\Zc G$.
\end{thm0} 

There is also a version for the automorphism group.

\begin{thm0}\label{aut ses}\cite{MalfaitWim2002Toag}
	There exist short exact sequences
	\[\begin{tikzcd}
		1\ar{r}{}&\mathcal{K}'\ar{r}{}&\Aut(G,K)\ar{r}{}&L'_1\ar{r}{}&1,
	\end{tikzcd}\]
and
	\[\begin{tikzcd}
		1\ar{r}{}&Z_\psi^1(Q,\Zc K)\ar{r}{}&\mathcal{K}'\ar{r}{}&L'_2\ar{r}{}&1
	\end{tikzcd}\]
	such that $L'_1\leq \Aut(Q)$ and $L'_2\leq \Aut(K)$.
\end{thm0} 

\begin{remark0}
If $K$ is a characteristic subgroup of $G$ then $\Aut(G,K)=\Aut(G)$ and $\Out(G,K)=\Out(G)$.
\end{remark0}

The relation between the Minkowski property and short exact sequences is explained in the following lemma:

\begin{lemma0}\label{Minkowski ses}
	Let $\begin{tikzcd}		1\ar{r}{}& K\ar{r}{}&G\ar{r}{p}&Q\ar{r}{}&1
	\end{tikzcd}$ be a short exact sequence of groups. If $K$ and $Q$ are Minkowski, then $G$ is Minkowski. If $K$ is finite and $G$ is Minkowski then $Q$ is Minkowski.
\end{lemma0}

\begin{proof}
Let us prove the first part of the statement. Let $F$ be a finite subgroup of $G$. Then we have a short exact sequence $1\longrightarrow F_1\longrightarrow F\longrightarrow F_3\longrightarrow 1$, where $F_3$ is the image of $F$ by the map $G\longrightarrow Q$ and $F_1=F\cap K$. If $C_1$ and $C_3$ are the Minkowski constants of $K$ and $Q$ respectively then $|F|\leq |F_1||F_3|\leq C_1C_3$. Therefore $G$ is Minkowski.

For the second part assume that $F$ is a finite subgroup of $Q$, then $p^{-1}(F)$ is a finite subgroup of $G$ of order $|K||F|$. Since $G$ is Minkowski we have that $|K||F|\leq C_2$ and therefore $|F|\leq C_2/|K|$. Thus $Q$ is Minkowski. 
\end{proof}

In particular, virtually torsion-free groups are Minkowski. A first consequence of these results is the following corollary:

\begin{corollary0}\label{outer automorphism and finite index subgroups}
	Let $\Gamma$ and $\Gamma'$ be finitely generated groups with finitely generated center such that $\Gamma'\trianglelefteq\Gamma$ and $\Gamma/\Gamma'=F$ is a finite group. If $\Out(\Gamma')$ is a Minkowski, then $\Out(\Gamma)$ is Minkowski.
\end{corollary0}

\begin{proof}
	By \cite[Lemma 1.(a)]{mccool1988automorphism} we know that $[\Out(\Gamma):\Out(\Gamma,\Gamma')]<\infty$, hence it is enough to prove that $\Out(\Gamma,\Gamma')$ is Minkowski.
	
	By \cref{out ses} and \cref{Minkowski ses}, the group $\Out(\Gamma,\Gamma')$ is Minkowski if $\Out(F)$, $C_{\Out(\Gamma')}\psi(F)/\psi(\Zc F)$ and $H^1_\psi(F,\Zc\Gamma')$ are Minkowski. But $\Out(F)$ and $H^1_\psi(F,\Zc\Gamma')$ are Minkowski since $F$ is finite and $\Zc\Gamma'$ is finitely generated. Lastly, $C_{\Out(\Gamma')}\psi(F)/\psi(\Zc F)$ is Minkowski since $C_{\Out(\Gamma')}\psi(F)\leq \Out(\Gamma')$ is Minkowski by hypothesis and $\psi(\Zc F)$ is a finite group, hence we can use the second part of \cref{Minkowski ses}. 
\end{proof}
\section{Finite group actions on aspherical manifolds: proof \cref{main theorem1}}\label{sec aspherical manifolds}

Let $F$ be a finite group acting effectively on a connected manifold $M$. Recall that if the action of $F$ has a fix point $x\in M$, then the group action induces a group morphism $F\longrightarrow \Aut(\pi_1(M,x))$. However, in the general case were the action does not have a fix point this group morphism is only well defined up to conjugation, so we have instead a group morphism $\psi:F\longrightarrow \Out(\pi_1(M))$. This group morphism is specially relevant when $M$ is a closed aspherical manifold, as the following theorem of P.E. Conner and F. Raymond shows:

\begin{thm0}\cite[Theorem 3.1.16]{lee2010seifert}\label{K-manifolds inner abelain}
Let $M$ be a closed aspherical manifold such that $\Zc(\pi_1(M))$ is finitely generated. Assume that $F$ is a finite group acting effectively on $M$. Then:
\begin{itemize}
	\item[1.] $\ker\psi$ is isomorphic to a subgroup of the torus $T^k$, where $k$ is the rank of $\Zc(\pi_1(M))$.
	\item[2.] If $F$ has a fixed point $x$ then the group morphism $F\longrightarrow \Aut(\pi_1(M,x))$ is injective.
\end{itemize}
  
\end{thm0}

This theorem was used to find the first examples of asymmetric manifolds, by constructing a closed aspherical manifold with a fundamental group with trivial center and torsion-free outer automorphisms group. In our case, it is the main tool to prove \cref{main theorem1}.  

\begin{proof}[Proof of part 1. of \cref{main theorem1}]
Let $C$ denote the bound of the order of finite subgroups of $\Out(\pi_1(M))$. If $F$ is a finite group acting effectively on $M$, then $\ker\psi$ is an abelian subgroup of $F$ and $[F:\ker\psi]=|F/\ker\psi|\leq C$ since $F/\ker\psi\leq\Out(\pi_1(M))$. Thus $\Homeo(M)$ is Jordan.
\end{proof}

To prove the second part we need the following group theoretic results (the second one due to Schur).

\begin{lemma0}\label{group theory discsym}\cite[Lemma 2.1]{riera2021topological}
Let $a,b,C$ be natural numbers and suppose that $F$ is a subgroup of $(\Z/a)^b$ such that $[(\Z/a)^b:F]\leq C$. Then there exists a subgroup $F'\leq F$ isomorphic to $(\Z/a')^b$ such that $C!a'\geq a$.
\end{lemma0}

\begin{lemma0}\label{schur inner theorem}\cite[Theorem 4.12]{robinson2013finiteness}
Let $\Gamma$ be a finitely generated group such that $[\Gamma:\Zc\Gamma]<\infty$. Then the commutator subgroup $[\Gamma,\Gamma]$ is finite.
\end{lemma0}

\begin{proof}[Proof of part 2. of \cref{main theorem1}]
We start proving the inequality $\D(M)\leq k$. Let $\{a_i\}_{i\in\mathbb{N}}$ be a sequence of natural numbers such that $a_i\longrightarrow\infty$ and $(\Z/a_i)^b$ acts effectively on $M$ for some $b\in\mathbb{N}$. We have induced group morphisms $\psi_i:(\Z/a_i)^b\longrightarrow \Out(\pi_1(M))$ for each $i$ such that $[(\Z/a_i)^b:\ker\psi_i]\leq C$.  By \cref{group theory discsym}, there exists a sequence $\{a'_i\}_{i\in\mathbb{N}}$ such that $(\Z/a'_i)^b\leq\ker\psi_i$. Since $a_i'C!\geq a_i$ we have that $a'_i\longrightarrow\infty$. Moreover, $(\Z/a'_i)^b$ is a subgroup of $T^k$ for any $i$, so we can conclude that $b\leq k$. In consequence $\D(M)\leq k$. To prove the inequality $k\leq n$ we take the $n$-dimensional manifold $\tilde{M}/\Zc(\pi_1(M))$, where $\tilde{M}$ is the universal cover of $M$. Since $\tilde{M}$ is contractible we have that $H^*(\tilde{M}/\Zc(\pi_1(M)))\cong H^*(T^k,\Z)$. The fact that $\tilde{M}/\Zc(\pi_1(M))$ has dimension $n$ implies that $H^i(\tilde{M}/\Zc(\pi_1(M)),\Z)=0$ for $i>n$, hence $k\leq n$.

Finally, we prove that $\D(M)=n$ if and only if $M\cong T^n$. It is clear that $\D(T^n)=n$ since $\Out(\Z^n)=\Gl(n,\Z)$ is Minkowski. Conversely, assume that $\D(M)=n$. Since the top cohomology $H^n(\tilde{M}/\Zc(\pi_1(M)),\Z)$ is non-zero, we can conclude that $\tilde{M}/\Zc(\pi_1(M))$ is a closed connected manifold and that the map $\tilde{M}/\Zc(\pi_1(M))\longrightarrow M$ is a regular finite cover. In consequence, $[\pi_1(M):\Zc\pi_1(M)]<\infty$ and $[\pi_1(M),\pi_1(M)]$ is trivial, by \cref{schur inner theorem} and the fact that $\pi_1(M)$ is torsion-free. Thus $\pi_1(M)\cong \Z^n$ and $M\cong T^n$, because the Borel conjecture holds for $\Z^n$. 
\end{proof}

The third part is a direct consequence of these two results:

\begin{lemma0}\label{aspherical euler char}\cite[Theorem IV.1.]{gottlieb1965certain} Let $X$ be a topological space with the homotopy type of a compact connected aspherical CW-complex. If $\chi(X)\neq 0$, then $\rank\Zc\pi_1(X)=0$.
\end{lemma0}

\begin{lemma0}\label{discsym almost asymetric}\cite[Lemma 8.1]{riera2023actions}
For any closed connected manifold $M$ we have $\D(M)=0$ if and only if $M$ is almost asymmetric.
\end{lemma0}

\begin{proof}[Proof of part 3. of \cref{main theorem1}]
If $M$ is a closed connected aspherical manifold, it has the homotopy type of a compact connected aspherical CW-complex. Thus, if $\chi(M)\neq 0 $ then $\D(M)\leq\rank\Zc(\pi_1(M))=0$. In consequence $\D(M)=0$ and $M$ is almost asymmetric.
\end{proof}

Finally, to prove the last part we need the following general lemma:

\begin{lemma0}\label{stab lemma}
	Let $M$ be a closed manifold. Assume that $\Homeo(M)$ is Jordan with constant $C$ and that there exists a constant $C'$ such that every finite group $F$ acting effectively on $M$ with a fixed point fulfils that $|F|\leq C'$. Then there exists a constant $D$ such that every finite group $F$ acting on $M$ has a subgroup $H$ such that $[F:H]\leq D$ and $|\Stab(H,M)|\leq D$.
\end{lemma0}

\begin{proof}
Let $F$ be a finite group acting effectively on $M$. Since $\Homeo(M)$ is Jordan there exists an abelian subgroup $H$ such that $[F:H]\leq C$. We will see that we can bound $\Stab(H,M)$ by  a constant $C''$ for any finite abelian group acting effectively on $M$ and therefore $D=\max\{C,C''\}$. Recall that there exists a constant $r$ such that any finite abelian group $H$ which acts effectively on $M$ fulfils that $\rank(H)\leq r$.

We have an inclusion $\Stab(H,M)\subseteq\{L\leq H:|L|\leq C'\}$, thus we will bound the number of subgroups of $H$ of order at most $C'$ instead of the number of stabilizers of the action of $H$. Since $H$ is abelian, we have a decomposition $H=H_1\times\dots\times H_l$, where each $H_i$ is a $p_i$-Sylow subgroup, thus is of the form $H_i=\Z/p_i^{a_{i,1}}\times \dots \Z/p_i^{a_{i,r_i}}$, where $a_{i,1}\leq \dots \leq a_{i,r_i}$ and $r_i\leq r$ for all $i$. Then we have a bijective correspondence 
$$\{L\leq H\} \longleftrightarrow \prod_{i=1}^l\{L_i\leq H_i\}$$ 
which induces an inclusion 
$$\{L\leq H:|L|\leq C'\} \subseteq \prod_{i=1}^l\{L_i\leq H_i:|L|\leq C'\}.$$

For a prime $p$ we define $e(p)=\max\{e\in\mathbb{N}:p^{e}\leq C'\}$ and we denote the set of all primes such that $e(p)\neq 0$ by $\mathcal{P}$. Note that $\mathcal{P}$ is precisely the set of primes which are equal or smaller than $C'$ and therefore $|\mathcal{P}|\leq C'$. Moreover, we have that 
$$\{L_i\leq H_i:|L|\leq C'\}\subseteq S_i=\{L_i\leq H_i:\text{all elements of $L$ have order at most $e(p_i)$}\}$$ for each $i$. Note that if $p_i\notin \mathcal{P}$ then $|S_i|=1$. Therefore $$\{L\leq H:|L|\leq C'\} \subseteq \prod_{i=1,p_i\in \mathcal{P}}^lS_i.$$

For each $H_i$ we choose an inclusion $H_i\longrightarrow (\Z/p_i^{a_{i,r_i}})^r$, which induce an inclusion 
$$S_i\subseteq \{L\leq (\Z/p_i^{a_{i,r_i}})^r:\text{all elements of $L$ have order at most $e(p_i)$}\}.$$

Finally, we use that if $L$ is a subgroup of $(\Z/p_i^{a_{i,r_i}})^r$ where all elements are of order at most $e(p_i)$ then $L_i\leq (\Z/p_i^{e(p_i)})^r\leq (\Z/p_i^{a_{i,r_i}})^r$. In consequence, we have an inclusion 
$$\{L_i\leq (\Z/p_i^{a_{i,r_i}})^r:\text{all elements of $L$ have order at most $e(p_i)$}\}\subset\{L_i\leq (\Z/p_i^{e(p_i)})^r\}.$$

Since $p_i\leq C'$ and $e(p_i)\leq C'$ we have that $|\{L_i\leq (\Z/p_i^{e(p_i)})^r\}|\leq 2^{|(\Z/p_i^{e(p_i)})^r|}\leq 2^{(C'^{C'})^r}$, which does not depend on the prime $p_i$. Therefore, we have that $|S_i|\leq 2^{(C'^{C'})^r}$ if $p_i\in \mathcal{P}$ and $|S_i|=1$ if $p_i\notin \mathcal{P}$. By using that $|\mathcal{P}|\leq C'$ we can conclude that $|\Stab(H,M)|\leq|\{L\leq H:|L|\leq C'\}|\leq  \prod_{i=1,p_i\in \mathcal{P}}^l|S_i|\leq (2^{(C'^{C'})^r})^{C'} $, which completes the proof.

\end{proof}

\begin{proof}[Proof of part 4. of \cref{main theorem1}]
We have already seen that $\Homeo(M)$ is Jordan. Since a finite group $F$ acting on $M$ with a fix point is a subgroup of $\Aut(\pi_1(M))$  and $\Aut(\pi_1(M))$ is Minkowski, the second condition of \cref{stab lemma} is also fulfilled. Thus, we only need to use \cref{stab lemma} to finish the proof of the fourth part.
\end{proof}

This ends the proof of \cref{main theorem1}. There are some interesting and natural questions that we summarise in the next remarks.

\begin{remark0}\label{remark Out Minkowski}
(The hypothesis on the fundamental group) When are the two hypothesis on the fundamental group fulfilled? No closed aspherical manifold with $\Zc(\pi_1(M))$ not finitely generated is known (see \cite[Remark 3.1.19.]{lee2010seifert}).
	
It is an interesting problem to find which closed connected aspherical manifolds $M$ satisfy that $\Out(\pi_1(M))$ is Minkowski (see \cite[Remark 7.1]{golota2023finite}). The next examples show some case where $\Out(\pi_1(M))$ is Minkowski:
\begin{itemize}
	\item[1.] In this paper we prove that aspherical locally homogeneous spaces (or classical aspherical manifolds following the terminology on \cite{farrell1990classical}) fulfil that $\Out(\pi_1(M))$ is Minkowski.
	\item[2.]  Another source of closed connected aspherical manifolds is the strict hyperbolization processes (see \cite{charney1995strict} and references therein). Given a closed oriented manifold $M'$ of dimension $n\geq 3$ we can construct a closed oriented aspherical manifold $M$ and a non-zero degree map $f:M\longrightarrow M'$ such that $\pi_1(M)$ is a hyperbolic group. These groups fulfils that $\Out(\pi_1(M))$ is finite (see \cite[5.4.A]{gromov1987hyperbolic} or \cite{paulin1991outer}). Moreover $\Zc(\pi_1(M))=0$ and therefore $M'$ is almost-asymmetric. 
	\item[3.] If $M$ is a closed connected aspherical 3-dimensional manifold, then $\Out(\pi_1(M))$ is Minkowski (see \cite{kojima1984bounding}). This fact was used in \cite{zimmermann2014jordan} to prove that $\Diff(M)$ is Jordan when $M$ is a closed smooth 3-manifold.
	\item[4.] If $\Out(\pi_1(M))$ is a finitely generated virtually abelian group then $\Out(\pi_1(M))$ is Minkowski. This is the case for piecewise linear locally symmetric spaces, a new type of closed aspherical manifold defined in \cite{tam2011gluing} by compactifying non-compact symmetric spaces and using the reflection trick on its corners. If $M$ is a piecewise linear locally symmetric space then $\Out(\pi_1(M))$ is finitely generated and virtually abelian by \cite[Theorem 4]{tam2011gluing}. Moreover, $\Zc\pi_1(M)=\{e\}$ by \cite[Lemma 7]{tam2011gluing}, hence piecewise linear locally symmetric spaces are almost asymmetric.
\end{itemize}

\end{remark0}


\begin{remark0}\label{remark out not Minkowski}
(Groups whose (outer) automorphism group is not Minkowski) There exists groups $\Gamma$ such that $B\Gamma$ is a finite CW-complex and $\Out(\Gamma)$ is not Minkowski. For example, let $\Gamma$ be the Baumslag-Solitar group $B(m,ml)=\langle a,b| ba^mb^{-1}=a^{ml}\rangle$ with $m,l\geq 2$. The space $BB(m,ml)$ is a finite aspherical $2$-dimensional CW-complex since it is a torsion-free one relator group (see \cite{lyndon1977combinatorial}). Moreover $\Out(B(m,ml))$ and $\Aut(B(m,ml))$ have elements of order $l^t(l-1)$ for arbitrarily large $t$ (see \cite[Lemma 3.8]{collins1983automorphisms} or \cite{levitt2007automorphism}) and thus they are not Minkowski.
\end{remark0}

\begin{remark0}\label{remark discsym=center}
(The discrete degree of symmetry versus the toral degree of symmetry) The question of when $\A(M)$ is equal to $\rank\Zc(\pi_1(M))$ has been extensively studied for closed connected aspherical manifolds $M$. We also note that any effective torus action on a closed aspherical manifold $M$ is almost-free (see \cite[Corollary 3.1.12]{lee2010seifert}), therefore $\A(M)$ is equal to the toral rank of $M$ (see \cite[\S 7.3]{felix2008algebraic}). 

If $\Out(\pi_1(M))$ is Minkowski and $\A(M)=\rank\Zc\pi_1(M)$ then $\D(M)=\rank\Zc\pi_1(M)$, since $\A(M)\leq\D(M)\leq \rank \Zc(\pi_1(M))$. The equality holds for infra-solvmanifolds or some aspherical locally homogeneous spaces (see \cite[Section 11.7]{lee2010seifert}). In \cref{submain theorem2 discsym} we prove that the equality is valid for closed aspherical locally homogeneous spaces. On the other hand, there exist closed aspherical manifolds such that $\A(M)=0$ and $\D(M)\geq 1$ (see \cite{cappell2013closed,riera2021topological}). It is an interesting question whether all closed connected aspherical manifolds satisfy $\D(M)=\rank\Zc\pi_1(M)$.  
\end{remark0}

\begin{remark0}
(Euler characteristic and asymmetry) The converse of part 3 of \cref{main theorem1} is not true. There exist asymmetric flat manifolds (see \cite{szczepanski2012geometry}). The fundamental group of a closed flat manifold fulfils the hypothesis of \cref{main theorem1} and the manifold Euler characteristic is $0$ since they are finitely covered by a torus.
\end{remark0}

\begin{remark0}
(Generalizations of aspherical manifolds) \Cref{K-manifolds inner abelain} has been generalized to closed connected manifolds $M$ satisfying that the only periodic self-homeomorphisms of the universal cover $\tilde{M}$ commuting with the deck of transformation groups $\pi_1(M)$ are elements of $\Zc\pi_1(M)$ (following the terminology of \cite[Theorem 3.2.2]{lee2010seifert} we will say that $M$ is admissible). This class of manifolds includes, for example, manifolds which admits a non-zero degree map to a closed aspherical manifold.

If $M$ is a closed connected admissible manifold such that $\Out(\pi_1(M))$ is Minkowski and $\Zc\pi_1(M)$ is finitely generated then $\Homeo(M)$ is Jordan and $\D(M)\leq\rank(\Zc\pi_1(M)/\operatorname{Torsion}(\Zc\pi_1(M)))$. For example, if $M$ is a closed aspherical manifold satisfying the hypothesis of \cref{main theorem1} and $N$ is a closed simply-connected manifold of the same dimension as $M$, then $M\# N$ is a closed admissible manifold such that $\Homeo(M\#N)$ is Jordan and $\D(M\# N)\leq \rank\Zc\pi_1(M)$. This last inequality is usually strict. For example, assume that $N$ is an odd dimensional manifold such that $\chi(N)\neq 0$. Then $\chi(T^n\#N)\neq 0$ and hence $\D(T^n\#N)=0<n=\rank \Zc(\pi_1(T^n\#N))$ (see \cite{riera2021topological}).

However, there exists a closed admissible manifold $M$ such that $\Homeo(M)$ is Jordan but $\Out(\pi_1(M))$ is not Minkowski. In \cite[Theorem 1.1.1]{bloomberg1975manifolds} it is proved that if $\Gamma$ and $\Lambda$ are not free groups then $\Out(\Gamma*\Lambda)\cong\Aut(\Gamma)\times\Aut(\Lambda)$ if $\Gamma\ncong\Lambda$ and $\Out(\Gamma*\Lambda)\cong(\Aut(\Gamma)\times\Aut(\Lambda))\times\Z/2$ if $\Gamma\cong\Lambda$. In particular, if $M$ is a closed 4-manifold with $\pi_1(M)\cong B(m,ml)$ with $m,l\geq 2$ then $\Out(\pi_1(M\#T^4))\cong\Aut(B(m,ml))\times\Gl(4,Z)$ is not Minkowski (see \cref{remark Out Minkowski}), but $M\#T^4$ is hypertoral and hence $\Homeo(M\#T^4)$ is Jordan by \cite[Theorem 1.15]{riera2021topological}.
\end{remark0}

\begin{remark0}
(Non-compact aspherical manifolds) If $M$ is a non-compact connected aspherical manifold then $\Homeo(M)$ is not necessarily Jordan, even when $\Out(\pi_1(M))$ is Minkowski. For example, since $\R^3$ admits effective actions by $\SO(3)$ then $\Homeo(T^2\times \R^3)$ is not Jordan by \cite[Theorem 1]{riera2017non}.
\end{remark0}


We end this section by exposing some facts on the relation between large finite group actions on manifolds and covering maps.

\begin{lemma0}\label{Jordan finite covering}
Let $M$ and $M'$ be closed connected manifolds and $p:M'\longrightarrow M$ be a finite covering. Then: 
\begin{itemize}
	\item[1.] If $\Homeo(M')$ is Jordan, then $\Homeo(M)$ is Jordan. 
	\item[2.] $\D(M')\geq\D(M)$.
	\item[3.] Assume that there exists a constant $D'$ such that any finite group $F'$ acting effectively on $M'$ with a fix point satisfies $|F'|\leq D'$. Then there exists a constant $D$ such that any finite group $F$ acting effectively on $M$ with a fix point satisfies $|F|\leq D$
\end{itemize}
\end{lemma0}

\begin{proof}
The first two parts are proven in \cite[\S 2.3]{IgnasiMundetiRiera2010Jtft} and  \cite[Theorem 1.12]{riera2021topological}. The proof of the third part follows the same arguments as the proofs of the first two parts.

Assume that $p:M'\longrightarrow M$ is a $n$-sheeted covering and $F$ is a finite group acting effectively on $M$. Then $F$ also acts on $\operatorname{Cov}_n(M)$, the set of $n$-sheeted coverings of $M$, by pull-backs. On the other hand  $\operatorname{Cov}_n(M)\cong\Hom(\pi_1(M),S_n)/\sim$ where $S_n$ is the $n$-th symmetric group and the equivalence relation is given by conjugation of elements of $S_n$. Therefore $\operatorname{Cov}_n(M)$ is finite, which implies that there exists a constant $C$ only depending on $M$ and $n$ such that any finite group $F$ acting effectively on $M$ has a subgroup $F_0$ which acts trivially on $\operatorname{Cov}_n(M)$ and $[F:F_0]\leq C$. Then there exists a finite group $F'_0$ acting effectively on $M'$ and a surjective group morphism $\pi:F'_0\longrightarrow F_0$ which makes the covering map $p:M'\longrightarrow M$ $\pi$-equivariant and $|\ker\pi|\leq n!$. 

Let $F$ be a finite group acting effectively on $M$ with a fix point $x\in M$. Then $x$ is also fixed by the action of $F_0$ and $F'_0$ acts on $p^{-1}(x)$. Given $x'\in p^{-1}(x)$, we have that the orbit of $x'$ by $F'_0$ is $F_0'/{F'_0}_{x'}\subseteq p^{-1}(x)$, which implies that $|F_0'/{F'_0}_{x'}|\leq n$. Finally, ${F'_0}_{x'}$ acts effectively on $M'$ with a fixed point, therefore $|{F'_0}_{x'}|\leq D'$. Since $\pi$ is surjective, we can take $D=C\cdot D' \cdot n$.
\end{proof}

In consequence, if $p:M' \longrightarrow M$ is a covering of closed connected aspherical manifolds and $M'$ fulfils the hypothesis of \cref{main theorem1} then all the conclusions of \cref{main theorem1} also hold for $M$. This fact can also be deduced for regular coverings using \cref{outer automorphism and finite index subgroups}.

Given a regular covering $p:M'\longrightarrow M$, an interesting problem is to determine when $\D(M)=\D(M')$. The next proposition will be used to give a partial answer to this question for the case of closed connected aspherical manifolds, which will be given in \cref{discsym coverings corollary}.

\begin{prop0}\label{discsym coverings}
Let $M$ be a closed connected aspherical manifold such that $\Zc(\pi_1(M))$ is finitely generated. Assume that $F$ is a finite group acting freely on $M$. Then $k=\rank(\Zc(\pi_1(M)))=\rank(\Zc(\pi_1(M/F)))$ if and only if the map $\psi':F\longrightarrow \Out(\Zc(\pi_1(M)))$ is trivial.
\end{prop0}

\begin{proof}
Recall that the free action of $F$ on $M$ induces a commutative diagram
\[\begin{tikzcd}
	& 1\ar{d}{}& 1 \ar{d}{}& 1 \ar{d}{}& \\
1\ar{r}{}& \Zc\pi_1(M)\ar{d}{}\ar{r}{}& C_{\tilde{F}}(\pi_1(M)) \ar{d}{}\ar{r}{}& \ker\psi \ar{d}{}\ar{r}{}& 1\\
1\ar{r}{}& \pi_1(M)\ar{d}{}\ar{r}{}& \tilde{F} \ar{d}{\tilde{\psi}}\ar{r}{p}& F \ar{d}{\psi}\ar{r}{}& 1\\
1\ar{r}{}& \Inn\pi_1(M)\ar{d}{}\ar{r}{}& \tilde{\psi}(\tilde{F}) \ar{d}{}\ar{r}{}& \psi(F) \ar{d}{}\ar{r}{}& 1\\	
	& 1& 1 & 1 & \\
\end{tikzcd}\]
where $\tilde{\psi}(\tilde{F})\leq \Aut(\pi_1(M))$ and $\psi(F)\leq\Out(\pi_1(M))$. We also note that $\Zc\tilde{F}\trianglelefteq C_{\tilde{F}}(\pi_1(M))$ and $\tilde{F}=\pi_1(M/F)$ is torsion-free. Recall also that $\tilde{\psi}(\tilde{f})=c_{\tilde{f}|\pi_1(M)}$ and that ${\psi}(f)=[c_{\sigma{f}|\pi_1(M)}]$ where $\sigma: F\longrightarrow \tilde{F}$ is a set-theoretic section of $p$.

Since $\Zc\pi_1(M)$ is a characteristic subgroup of $\pi_1(M)$, we have maps $\tilde{\psi}':\tilde{F}\longrightarrow \Gl(k,\Z)$ and $\psi':F\longrightarrow \Gl(k,\Z)$ by restricting (outer) automorphisms to $\Zc\pi_1(M)$. Moreover, $\tilde{\psi}'=\psi'\circ p$ and since $p$ is surjective we can conclude that $\tilde{\psi}'$ is trivial if and only if $\psi'$ is trivial.

Assume now that $\psi'$ is trivial and hence $\tilde{\psi}'$ is trivial too. This implies that $[\tilde{f},z]=e$ for any $\tilde{f}\in\tilde{F}$ and any $z\in \Zc\pi_1(M)$. Thus $\Zc\pi_1(M)\trianglelefteq \Zc\tilde{F}\trianglelefteq C_{\tilde{F}}(\pi_1(M))$. But $\rank(\Zc\pi_1(M))=\rank(C_{\tilde{F}}(\pi_1(M)))$ since the first row of the commutative diagram is a central exact sequence and $C_{\tilde{F}}(\pi_1(M))$ is torsion-free. In consequence, $\rank(\Zc\pi_1(M))=\rank(\Zc{\tilde{F}})$.

If $\rank(\Zc\pi_1(M))=\rank(\Zc{\tilde{F}})$ then $[C_{\tilde{F}}(\pi_1(M)):\Zc\tilde{F}]< \infty$. This implies that $\tilde{\psi}'(\tilde{F})\leq \Gl(r,\Z)$ fixes a sublattice of $\Zc\pi_1(M)$ and therefore $\tilde{\psi}'$ is trivial. Thus $\psi'$ is also trivial, as desired.
\end{proof}

\begin{corollary0}\label{discsym coverings corollary}
Let $M$ be a closed connected aspherical manifold such that $\Zc\pi_1(M)$ is finitely generated, $\Out(\pi_1(M))$ is Minkowski and $\D(M)=\rank\Zc(\pi_1(M))$. Assume that $F$ is a finite group acting freely on $M$ such that $\psi'$ is not trivial. Then $\D(M/F)<\D(M)$.
\end{corollary0}

See \cref{remark discsym=center} for examples where the hypothesis of the corollary holds true. 

\section{Solvmanifolds}\label{sec:solvmanifolds}

Recall that a manifold $M$ is a solvmanifold if $M$ admits transitive group action of a connected solvable Lie group $R$. If $R$ is nilpotent then we will say that $M$ is a nilmanifold. The structure of nilmanifolds and solvmanifolds and their fundamental groups have been vastly studied (see \cite{maltsev1951class,GorbatsevichVV2009Csod,auslander1973} and \cite[Chapter II,III,IV]{raghunathan2012discrete}). The fundamental group of any solvmanifold is polycyclic. Recall that $\Gamma$ is a polycyclic group if there exists a sequence of $\Gamma\trianglerighteq\Gamma_1\trianglerighteq\dots\trianglerighteq\Gamma_r=\{e\}$ such that each quotient $\Gamma_{i-1}/\Gamma_i$ is cyclic. Equivalently, a group is polycyclic if it is solvable group and all its subgroups are finitely generated. In particular, the center of a polycyclic group is finitely generated. Recall that $\Gamma$ is a virtually polycyclic group if it contains a polycyclic subgroup of finite index. 

\begin{thm0}\cite{Wehrfritz1994}\label{outer automorphisms polycyclic} Let $\Gamma$ be a virtually polycyclic group, then $\Out(\Gamma)$ is isomorphic to a subgroup of $\Gl(n,\mathbb{Z})$ for some $n$.
\end{thm0}

The proof of this theorem uses the analogue statement for the automorphism group.

\begin{thm0}\cite{mebzljakov1970integral}\label{automorphisms polycyclic}
Let $\Gamma$ be a virtually polycyclic group, then $\Aut(\Gamma)$ is isomorphic to a subgroup of $\Gl(n,\mathbb{Z})$ for some $n$.
\end{thm0}

Hence, if $\Gamma$ is virtually polycyclic, then $\Out(\Gamma)$ and $\Aut(\Gamma)$ are Minkowski. In consequence, by using \cref{main theorem1} we can conclude that a closed solvmanifold $M$ has Jordan homeomorphism group, $\D(M)\leq \rank(\Zc(\pi_1(M)))$ (and indeed $\D(M)= \rank(\Zc(\pi_1(M)))$ by the results in \cite[Section 11.7]{lee2010seifert}) and there exists a constant $C$ such that any finite group $F$ acting effectively on $M$ has a subgroup $H$ such that $|\Stab(H,M)|\leq C$ and $[F:H]\leq C$. These properties are also satisfied by any manifold finitely covered by a solvmanifold. This includes flat manifolds, almost-flat manifolds and infra-solvmanifolds.

\begin{example0}
Let us show a low dimensional example in order to illustrate some consequences of \cref{main theorem1} on nilmanifolds. The 3-dimensional Heisenberg group 
\[H=\{(x,y,z)=\left(\begin{array}{ccc} 1 & x & z\\  0 & 1 & y \\  0 & 0 & 1 \end{array} \right):x,y,z\in \R\} \] is a simply connected nilpotent Lie group. Any lattice of $H$ is isomorphic to a lattice of the form
\[\Gamma_k=\{(x,y,z)=\left(\begin{array}{ccc} 1 & x & \tfrac{1}{k}z\\  0 & 1 & y \\  0 & 0 & 1 \end{array} \right):x,y,z\in \Z\} \]
where $k$ is a positive integer. Note that $\Zc\Gamma_k=\langle(0,0,1)\rangle$ and two lattices $\Gamma_k\cong\Gamma_l$ are isomorphic if and only if $k=l$. A possible presentation of these lattices is $\Gamma_k=\langle a,b,c| [c,a]=[c,b]=1, [a,b]=c^k \rangle$ where $a=(1,0,0)$, $b=(0,1,0)$ and $c=(0,0,1)$. Thus $\Zc\Gamma_k=\langle c\rangle\cong\Z$ for all $k$. In \cite[\S 8]{conner2006manifolds} it is shown that $\Aut(\Gamma_k)\cong \Z^2\rtimes\Gl(2,\Z)$ and $\Out(\Gamma_k)\cong (\Z/k)^2\rtimes\Gl(2,\Z)$. Note that as an abstract group, $\Aut(\Gamma_k)$ does not depend on $k$, but $\Out(\Gamma_k)$ does. Both of them are Minkowski and therefore the conclusions of \cref{main theorem1} are valid for $H/\Gamma_k$. In particular, $\Homeo(H/\Gamma_k)$ is Jordan and $\D(H/\Gamma_k)=1$. However, not all finite subgroups of $\Out(\Gamma_k)$ can be realized by a group action on $H/\Gamma_k$ (see \cite{raymond1977failure}).

If a $n$-dimensional flat manifold such that $\D(M)=n=\D(T^n)$ then $M$ is homeomorphic to $T^n$. On the other hand, an almost-flat manifold $M$ finitely covered by a nilmanifold $N/\Gamma$ satisfying that $\D(M)=\D(N/\Gamma)$ is not necessarily isomorphic to $N/\Gamma$. The 3-dimensional Heisenberg manifold can be used to construct an almost-flat manifold $M$ which is finitely covered by $H/\Gamma_2$ and $\D(M)=\D(H/\Gamma_2)=1$ but $M$ is not homeomorphic to $H/\Gamma_2$ ($M$ is not even a nilmanifold). There is a free action of $\Z/2$ on $H/\Gamma_2$ such that its orbit space $M$ is an almost-flat manifold with fundamental group $\pi_1(M)=\langle a,b,c,\alpha| [c,\alpha]=[c,a]=[c,b]=1, [a,b]=c^2,\alpha a=a^{-1}\alpha, \alpha b=b^{-1}\alpha,\alpha^2=c \rangle$ (see \cite[pg. 160]{dekimpe2006almost}). It is clear that $\langle c\rangle\leq\Zc\pi_1(M)$ and therefore $1\leq \D(M)\leq \D(H/\Gamma_2)=1$. In consequence, we have an equality $\D(M)= \D(H/\Gamma_2)=1$. Note that $\alpha c\alpha^{-1}=c$, hence the morphism $\psi':\Z/2\longrightarrow \Out(\Zc\Gamma_2)$ is trivial, as we expected from \cref{discsym coverings corollary}. 

Let $M$ be a closed aspherical manifold satisfying the hypothesis of \cref{main theorem1}. The fact that $\D(M)=\rank(\Zc(\pi_1(M)))$ does not imply that there are no effective actions of $(\Z/k)^r$ with $r>\D(M)$ for some $k$. Indeed, if $M$ is a compact solvmanifold $R/\Gamma$ then \cite[Corollary 3.3]{JoLee} asserts that if $(\Z/p)^r$ acts freely on $R/\Gamma$ then $r\leq \dim R/\Gamma$. This bound is sharp even when $R/\Gamma$ is not a torus. For an integer $k\leq 2$ we consider the subgroup 

\[\Gamma'_k=\{(x,y,z)=\left(\begin{array}{ccc} 1 & \tfrac{1}{k}x & \tfrac{1}{k^3}z\\  0 & 1 & \tfrac{1}{k}y \\  0 & 0 & 1 \end{array} \right):x,y,z\in \Z\} \]	
which is a lattice of $H$ isomorphic to $\Gamma_k$. An straightforward computation shows that $\Gamma_{k^2}$ is a normal subgroup of $\Gamma'_k$ and $\Gamma'_k/\Gamma_{k^2}\cong(\Z/k)^3$. Therefore $(\Z/k)^3$ acts freely on $H/\Gamma_{k^2}$ even though $\D(H/\Gamma_{k^2})=1$ (see \cite{choi2005free} for a classification of all finite abelian group actions on Heisenberg manifolds). Therefore the bound of \cite[Corollary 3.3]{JoLee} is sharp. 
\end{example0}

\section{Locally symmetric spaces}\label{sec:locally symmetric spaces}

A locally symmetric space is a double coset space $H\setminus G/\Gamma$ where $G$ is a connected semisimple Lie group, $H$ is a maximal compact subgroup of $G$ and $\Gamma$ is a torsion-free lattice of $G$. Since $H\setminus G\cong \R^n$ for some $n$ (see \cite[Chapter VI, Theorem 1]{helgason2001differential}), then $H\setminus G/\Gamma$ is a connected aspherical manifold with fundamental group $\Gamma$. If we assume that $H\setminus G/\Gamma$ is compact then $G/\Gamma$ is also compact and therefore $\Gamma$ is a cocompact lattice. 

The aim of this section is to prove the following proposition.

\begin{prop0}\label{out lattice finite}
Let $\Gamma$ be a lattice of a connected semisimple Lie group $G$ without compact factors. Then $\Out(\Gamma)$ Minkowski.
\end{prop0}

The objective of the first part of the section is to recall the results of the theory of lattices of semisimple Lie groups and the theory of relatively hyperbolic groups needed to prove \cref{out lattice finite}.  
We start stating some well-known facts and theorems on lattices of semisimple Lie groups and hyperbolic groups. We refer to \cite{morris2001introduction} for an introduction to this topic. 

\begin{thm0}\label{facts semisimple lattices} Let $G$ be a connected semisimple linear group without compact factors and $\Gamma$ a lattice of $G$. Then:
	\begin{itemize}
		\item[1.] (Selberg's lemma, \cite[(4.2.8) Theorem]{morris2001introduction}) Any lattice $\Gamma\leq G$ contains a normal torsion-free lattice $\Gamma'$ such that $[\Gamma:\Gamma']<\infty$.
		\item[2.] (Reducible lattices, \cite[(4.3.3) Proposition]{morris2001introduction}) Assume that $G$ is centreless. Then there exist simple subgroups $G_1,...,G_r$ of $G$ and lattices $\Gamma_i\leq G_i$ for all $1\leq i\leq r$ such that $G=G_1\times \cdots G_r$ and $\Gamma_1\times\cdots\times \Gamma_r$ is a normal finite index subgroup of $\Gamma$.
		\item[3.] (Consequences of Borel density theorem I, \cite[(4.5.3) Corollary]{morris2001introduction}) The group $N_G(\Gamma)/\Gamma$ is finite and $C_G(\Gamma)=\Zc G$ (equivalently, $\Zc\Gamma=\Gamma\cap\Zc G$).
		\item[4.] (Consequences of Borel density theorem II, \cite[(4.5.2) Corollary]{morris2001introduction}) Let $H$ be a connected closed subgroup of $G$. Assume that $\Gamma$ normalizes $H$, then $H$ is normal in $G$.
		\item[5.] (Mostow-Prasad-Margulis rigidity theorem, \cite[(15.1.2) Thereom]{morris2001introduction}) Let $G'$ be another connected semisimple linear Lie group without compact factors, and let $\Gamma'$ be a lattice of $G'$. Assume that $G$ and $G'$ have trivial center. Finally, assume that there does not exists any simple factor $H$ of $G$ such that $H\cong \operatorname{PSL}(2,\R)$ and $H\cap\Gamma_1$ is a lattice in $H$. Then any isomorphism from $\Gamma$ to $\Gamma'$ extends to a unique continuous isomorphism from $G$ to $G'$.
		\item[6.] (Margulis superrigidity, \cite[\S 16.1]{morris2001introduction}) Assume that $\rank_\R G\geq 2$ and that $\Gamma$ is an irreducible lattice. Given a representation $\rho:\Gamma\longrightarrow \Gl(n,\R)$, let $\overline{\rho(\Gamma)}^o$ be the identity component of the Zariski closure of $\rho(\Gamma)$ in $\Gl(n,\R)$ and set $\Gamma_0=\rho^{-1}(\overline{\rho(\Gamma)}^o)$. Then there exists a representation $\tilde{\rho}:G\longrightarrow \overline{\rho(\Gamma)}^o$ such that $\tilde{\rho}_{|\Gamma_0}=\rho_{|\Gamma_0}$. 
		\item[7.] (Margulis normal theorem, \cite[(17.1.1) Theorem]{morris2001introduction})	Assume that $G$ has $\rank_{\R}G\geq 2$ and finite center, and that $\Gamma$ is an irreducible lattice. If $N$ is a normal subgroup of $\Gamma$, then either $N\leq \Zc(G)$ or $\Gamma/N$ is finite.
\end{itemize}
\end{thm0}

From \cref{facts semisimple lattices}(3) we can deduce:

\begin{corollary0}\label{normal solvable groups lattices}
Let $G$ be a connected semisimple linear group without compact factors and $\Gamma$ a torsion-free lattice of $G$. Then $\Gamma$ does not contain normal solvable subgroups.  
\end{corollary0}

\begin{proof}
Assume on the contrary, that $\Gamma$ contains a solvable normal subgroup $\Lambda$. Then, the connected component of the identity of the Zariski closure $\overline{\Lambda}^o$ is a connected solvable Lie subgroup of $G$. Since $\Lambda$ is normal in $\Gamma$ then $\overline{\Lambda}^o$ is normalized by $\Gamma$. By \cref{facts semisimple lattices}(3), $\overline{\Lambda}^o$ is a connected solvable normal subgroup of $G$, which contradicts the semisimplicity of $G$.
\end{proof}

We now introduce the results we need from the theory of relatively hyperbolic groups. For an introduction to relatively hyperbolic groups we refer to \cite{bowditch1997relatively,osin2006relatively}. 

Let us recall the definition of relative hyperbolicity from \cite{osin2006relatively,minasyan2010normal}. Given a group $H$ and a collection of proper subgroups $\{H_i\}_{i\in I}$, a subset $X$ of $G$ is a relative generating set of $H$ with respect to $\{H_i\}_{i\in I}$ if $X$ together with the union of all $H_i$ generates $H$. In this situation $H$ can be written as a quotient group of a free group $\mathcal{F}=(*_{i\in I}H_i)*F(X)$, where $F(X)$ denotes the free group with basis $X$. If the kernel $\mathcal{F}\longrightarrow H$ is the normal closure of a subset $R$ of $\mathcal{F}$ then we say that $H$ has a relative presentation $\langle X,H_i,i\in I|R\rangle$. If $X$ and $R$ are finite sets then we say that $H$ is finitely presented relative to $\{H_i\}_{i\in I}$.

Set $\mathcal{H}=\bigsqcup_{i\in I}(H_i\setminus\{e\})$. Given a word $W$ in the alphabet $X^{\pm}\cup\mathcal{H}$ representing the trivial element $e$ in $H$ there exists an element of the form $\prod_{j=1}^{k}f_jR_j^{\pm 1}f_j^{-1}$ which is equal to $W$ in the group $\mathcal{F}$, where $R_j\in R$ and $f_j\in \mathcal{F}$. The smallest possible $k$ for which $W$ is equal to an element of the form $\prod_{j=1}^{k}f_jR_j^{\pm 1}f_j^{-1}$ is called the relative area of $W$ and denoted by $Area^{rel}(W)$. If $||W||$ denotes the length of $W$ in the alphabet $X^{\pm}\cup\mathcal{H}$ then:

\begin{defn0}\cite[Definition 2.1]{minasyan2010normal}\label{def relatively hyperbolic group}
	A group $H$ is hyperbolic relative to a collection of proper subgroup $\{H_i\}_{i\in I}$ if $H$ is finitely presented relative to $\{H_i\}_{i\in I}$ and there is a constant $C>0$ such that any word $W$ in $X^{\pm}\cup\mathcal{H}$ representing the identity in $H$ satisfies that $$Area^{rel}(W)\leq  C||W||.$$ 
	
	The groups $H_i$ are called peripheral subgroups. We say that a group $H$ is relatively hyperbolic if there exists a collection of subgroups $\{H_i\}_{i\in I}$ such that $H$ is hyperbolic relative to $\{H_i\}_{i\in I}$.
\end{defn0} 

\begin{remark0}
	The definition is independent of the choice of $X$ and $R$ (see \cite{osin2006relatively}).
\end{remark0}

\begin{remark0}
	There are other definitions of relative hyperbolicity of a group $H$ with respect a collection of subgroups $\{H_i\}_{i\in I}$ in the literature. If $H$ is torsion-free and finitely presented and $H_i$ is finitely presented for all $i$, then all the definitions are equivalent (see \cite[Definition 1.1]{belegradek2008endomorphisms} and references therein). This will happen in our setting.
\end{remark0}

\begin{remark0}
	A group $H$ is hyperbolic if it is hyperbolic relative to the collection of subgroups which only contains the trivial subgroup.
\end{remark0}

We need to introduce relative hyperbolicity because of the next theorem:

\begin{thm0}\label{lattices relatively hyperbolic}(\cite{farb1998relatively} and \cite[\S 0.2(F)]{gromov1987hyperbolic})
	A lattice $\Gamma$ in a connected semisimple Lie group without compact factors $G$ is relatively hyperbolic if and only if $\rank_\R G=1$. More precisely, $\Gamma$ is hyperbolic relative to the collection of all its cusp subgroups associated to the cusps of the symmetric space $H\setminus G/\Gamma$. The cusp subgroups are virtually nilpotent.
\end{thm0}

Let $H$ be a relatively hyperbolic group to $\{H_i\}_{i\in I}$. Recall that $h\in H$ is said to be parabolic if its is conjugate to an element of $H_i$ for some $i$. An element $h\in H$ which is not parabolic is said to be hyperbolic.

\begin{lemma0}\label{realtive hyperbolic elements}
	Let $H$ be torsion-free a relative hyperbolic group to $\{H_i\}_{i\in I}$. Then:
	\begin{itemize}
		\item[1.] \cite[Lemma 2.4]{minasyan2010normal} If $h\in H$ is hyperbolic then $C_H(h)$ is cyclic.
		\item[2.] \cite[Proposition 3.3]{minasyan2010normal} The set of hyperbolic elements generates $H$.
		\item[3.] \cite[Lemma 2.2]{minasyan2010normal} Given $i\in I$ and $h\in H\setminus H_i$, we have that $H_i\cap hH_ih^{-1}=\{e\}$.
	\end{itemize}
\end{lemma0}

We are interested in the following corollary of item 3.

\begin{corollary0}\label{centrailizer parabolic elements}
	Let $h$ be a non-trivial element of $H_i$, then $C_H(h)\leq H_i$.
\end{corollary0} 

\begin{proof}
	Suppose that $h'\in C_H(h)$. Then $h'hh'^{-1}=h$, therefore $H_i\cap h'H_ih'^{-1}\neq \{e\}$. By \cref{realtive hyperbolic elements}, $h'\in H_i$. Thus, $C_H(h)\leq H_i$. 
\end{proof}

We are ready to prove \cref{out lattice finite}. Firstly, let us state the following well-known result and give a proof of it for the sake of completeness.

\begin{lemma0}\label{out finite mostow}
Let $G$ be connected linear semisimple Lie groups with trivial center and no compact factors, let $\Gamma$ be a lattice of $G$ and assume that there does not exist any simple factor $H$ of $G$ such that $H\cong \operatorname{PSL}(2,\R)$ and $H\cap\Gamma_1$ is a lattice in $H$. Then the group $\Out(\Gamma)$ is finite. 
\end{lemma0}

\begin{proof}
Let $F:\Aut(\Gamma)\longrightarrow \Aut(G)$ be the morphism sending an automorphism of $\Gamma$ to its unique extension on $G$ by \cref{facts semisimple lattices}(5). Clearly, $F(\Inn(\Gamma))\leq \Inn(G)$, so $F$ descends to a group morphism $f:\Out(\Gamma)\longrightarrow \Out(G)$. Then, the claim follows from the fact that $\Out(G)$ is finite and that $\ker f=N_{G}(\Gamma)/\Gamma$ is also finite by \cref{facts semisimple lattices}(4).
\end{proof}

\begin{lemma0}\label{out Minkowski with center}
Let $G$ be a connected semisimple Lie groups without compact factors and let $\Gamma$ be a lattice of $G$. If $\Out(\Gamma/\Zc\Gamma)$ is Minkowski then $\Out(\Gamma)$ is Minkowski.
\end{lemma0} 

\begin{proof}
We consider the central short exact sequence $1\longrightarrow \Zc\Gamma\longrightarrow \Gamma\longrightarrow \Gamma/\Zc\Gamma\longrightarrow 1$. The center is a characteristic subgroup, hence $\Out(\Gamma)=\Out(\Gamma,\Zc\Gamma)$. 

On the other hand, $\Out(\Gamma/\Zc\Gamma)$ is Minkowski by hypothesis, $\Out(\Zc\Gamma)$ is Minkowski since $\Zc\Gamma$ is a finitely generated abelian group and $\overline{H}^1(\Gamma/\Zc\Gamma,\Zc\Gamma)$ is Minkowski since it is finitely generated and abelian. Therefore by \cref{out ses} we obtain that $\Out(\Gamma)$ is Minkowski. 
\end{proof}

In view of the preceding lemma, to prove \cref{out lattice finite} it remains to show that $\Out(\Gamma/\Zc\Gamma)$ is Minkowski. Note that $\Gamma/\Zc\Gamma$ is a lattice in $G/\Zc G$ which is a centreless connected semisimple linear Lie group without compact factors. By \cref{facts semisimple lattices} (1) and (2), $G/\Zc G=G_1\times \cdots \times G_n$ and $\Gamma/\Zc\Gamma$ has a normal torsion-free finite index subgroup $\Lambda$ of the form $\Lambda=\Lambda_1\times\cdots\times \Lambda_n$ where each $\Lambda_i$ is an irreducible lattice of $G_i$. Note that each of this groups is centreless. We will use the decomposition $\Lambda=\Lambda^{(1)}\times \Lambda^{(\geq 2)}$ where $\Lambda^{(1)}$ is the product of all lattices of real rank $1$ and $\Lambda^{(\geq 2)}$ is the product of all lattices of real rank greater or equal than $2$. After reordering we can assume that $\Lambda^{(1)}=\Lambda_1\times\cdots \Lambda_m$ and $\Lambda^{(\geq 2)}=\Lambda_{m+1}\times\cdots \Lambda_n$.

\begin{lemma0}\label{lattice rank 2 charactherisitc}
The group $\Lambda^{(\geq 2)}$ is characteristic in $\Lambda$.
\end{lemma0}

\begin{proof}
	We will prove that given any $f\in\Aut(\Lambda)$ and any $m+1\leq j\leq n$ we have $f(\Lambda_j)\leq \Lambda^{(\geq 2)}$. 
	
	Let $\pi_i:\Lambda\longrightarrow \Lambda_i$ and $\iota_i:\Lambda_i\longrightarrow \Lambda$ denote the natural projection and inclusion morphisms. We take the group morphism $\pi_i\circ f\circ \iota_j:\Lambda_j\longrightarrow \Lambda_i$ with $1\leq i\leq m$ and $m+1\leq j\leq n$. By Margulis normal subgroup theorem (see \cref{facts semisimple lattices}(7)) we have that $N=\ker(\pi_i\circ f\circ \iota_j)\trianglelefteq \Lambda_j$ is either finite or has finite index in $\Lambda_j$. If $N$ is finite then it is trivial since $\Lambda_j$ is torsion-free and $\pi_i\circ f\circ \iota_j$ is injective. If it has finite index then $\Lambda_j/N$ is a finite subgroup of $\Lambda_i$ and therefore it is trivial since $\Lambda_i$ is torsion-free. In this case, $\pi_i\circ f\circ \iota_j$ is trivial. 
	
	The morphisms $\pi_i\circ f\circ \iota_j:\Lambda_j\longrightarrow\Lambda_i$ cannot be injective for any $i$ and $j$. Indeed, since $G_i$ is centreless (hence $G_i\leq \Gl(n_i,\R)$ for some $i$), we can construct a representation $\rho:\Lambda_j\longrightarrow \Gl(n_i,\R)$ given by the composition $\pi_i\circ f\circ \iota_j$ with the inclusions $\Lambda_i\leq G_i\leq \Gl(n_i,\R)$. By Margulis superrigidity (see \cref{facts semisimple lattices}(6)) there exists a finite index subgroup $\Lambda_{j0}\leq \Lambda_j$ and a representation $\tilde{\rho}:G_j\longrightarrow \overline{\rho(\Lambda_j)}^o$ such that $\tilde{\rho}_{|\Lambda_{j0}}={\rho}_{|\Lambda_{j0}}$. Since  $\overline{\rho(\Lambda_j)}^o\leq G_i$, $\tilde{\rho} $ induces an group morphism $G_j\longrightarrow G_i$, which is trivial since $G_i$ and $G_j$ are simple, $\rank_R G_j\geq 2$ and $\rank_R G_i=1$. Therefore $\Lambda_{j0}\leq\ker\pi_i\circ f\circ \iota_j$. Since $\pi_i\circ f\circ \iota_j$ is not injective then $\pi_i\circ f\circ \iota_j$ is trivial, as desired.
\end{proof}

Note that $\Out(\Lambda^{(\geq 2)})$ is finite by \cref{out finite mostow}. 

\begin{lemma0}\label{out hyperbolic finite}
The group $\Out(\Lambda^{(1)})$ is Minkowski.
\end{lemma0}

\begin{proof}
	Let $H=\{f\in \Aut(\Lambda^{(1)}):f(\Lambda_i)=\Lambda_i\text{ for all $i$}\}$. Firstly, we will show that $H$ has finite index in  $\Aut(\Lambda^{(1)})$. Since $\Inn(\Lambda^{(1)})\trianglelefteq H$, this will imply that $[\Out(\Lambda^{(1)}):H/\Inn(\Lambda^{(1)})]<\infty$. 
	
	We denote by $r(\lambda)$ the number of non-trivial entries of an element $\lambda\in \Lambda^{(1)}$. Let $S$ be the set of elements of $\Lambda^{(1)}$ whose centralizer $C_{\Lambda^{(1)}}(\lambda)$ is isomorphic as an abstract group to $\Z\times \Lambda_\lambda$, where $\Lambda_\lambda$ is a product of lattices in $\R$-rank one semisimple Lie groups (which depend on the element $\lambda\in\Lambda$). 
	
	Firstly, we note that if $\lambda\in S$ then $r(\lambda)=1$. Given $\lambda=(\lambda_1,\dots,\lambda_m)\in \Lambda^{(1)}$, its centralizer is $C_{\Lambda^{(1)}}(\lambda)=\prod_{i=1}^{m}C_{\Lambda_i}(\lambda_i)$. In addition, $C_{\Lambda_i}(\lambda_i)=\Lambda_i$ if and only if $\lambda_i=e_i$. If $\lambda_i$ is not trivial then $C_{\Lambda_i}(\lambda_i)$ is virtually nilpotent. Indeed, if $\lambda_i$ is hyperbolic then $C_{\Lambda_i}(\lambda_i)\cong \Z$ by \cref{realtive hyperbolic elements}(1) and if $\lambda_i$ is parabolic then $C_{\Lambda_i}(\lambda_i)$ is virtually nilpotent by \cref{centrailizer parabolic elements} and \cref{lattices relatively hyperbolic}. If $r(\lambda)\neq 1$ and $\lambda\in S$ then $C_{\Lambda^{(1)}}(\lambda)$ would contain at least two virtually nilpotent factors and one of this factors would be a normal subgroup of a lattice in a centreless semisimple Lie group. This contradicts \cref{normal solvable groups lattices}.	Thus, we can conclude that $\lambda$ is not in $S$. However, note that it is possible for an element $\lambda$ with $r(\lambda)=1$ to not be in $S$. If $\lambda=(e_1,\dots,\lambda_i,\dots, e_m)$, then $C_{\Lambda^{(1)}}(\lambda)\cong C_{\Lambda_{i}}(\lambda_i)\times \Lambda_{\lambda}$ and $C_{\Lambda_{i}}(\lambda_i)$ is not necessarily isomorphic to $\Z$ if $\lambda_i$ is not hyperbolic.
	
	Clearly $S$ is preserved by automorphisms of $\Lambda^{(1)}$.  Finally, $S$ generates $\Lambda^{(1)}$ by \cref{realtive hyperbolic elements}(2), since $S$ contains all elements whose only non-trivial entry is hyperbolic. Consequently, any $f\in\Aut(\Lambda^{(1)})$ permutes the factors of $\Lambda^{(1)}$ and we can construct a group morphism $\phi:\Aut(\Lambda^{(1)})\longrightarrow S_m$ such that $H=\ker \phi$. Consequently, $[\Aut(\Lambda^{(1)}):H]<m!$.
	
	We are ready to prove that $\Out(\Lambda^{(1)})$ is Minkowski. We proceed by induction on the number of factors $m$. If $m=1$ then $\Lambda^{(1)}$ is an irreducible lattice in a centreless semisimple Lie group $G_1$. If $G_1\ncong \PSL(2,\R)$ then we can use \cref{out finite mostow} to conclude that $\Out(\Lambda^{(1)})$ is finite. If $G_1\cong \PSL(2,\R)$ then $\Lambda^{(1)}$ is a Fuchsian group and therefore $\Out(\Lambda^{(1)})$ is virtually torsion-free \cite[Corollary 2.6]{metaftsis2006residual} and hence Minkowski.
	
	Assume now that $\Out(\Lambda_1\times\cdots\times\Lambda_{m-1})$ is Minkowski. By \cref{out ses}, $\Out(\Lambda^{(1)},\Lambda_m)$ is Minkowski if $\Out(\Lambda_m)$, $\Out(\Lambda_1\times\cdots\times\Lambda_{m-1})$ and $H^1(\Lambda_1\times\cdots\times\Lambda_{m-1},\Zc\Lambda_m)$ are Minkowski. But they are Minkowski by induction hypothesis and the fact that $\Zc\Lambda_m$ is trivial. Since $H\leq \Aut(\Lambda^{(1)},\Lambda_m)\leq \Aut(\Lambda^{(1)})$ and $[\Aut(\Lambda^{(1)}):H]<\infty$ we have that $[\Out(\Lambda^{(1)}):\Out(\Lambda^{(1)},\Lambda_m)]<\infty$, which implies that $\Out(\Lambda^{(1)})$ is Minkowski.
\end{proof}

We are ready to prove that $\Out(\Gamma/\Zc\Gamma)$ is Minkowski and finish the proof of \cref{out lattice finite}. 

\begin{lemma0}\label{out cocompact minkowski}
The group $\Out(\Gamma/\Zc\Gamma)$ is Minkowski. 
\end{lemma0}

\begin{proof}
Since $[\Gamma/\Zc\Gamma:\Lambda]<\infty$, by \cref{outer automorphism and finite index subgroups} it is enough to prove that $\Out(\Lambda)$ is Minkowski. We have seen that  $\Out(\Lambda^{(\geq 2)})$, $\Out(\Lambda^{(1)})$ and $H^1(\Lambda^{(1)},\Zc\Lambda^{(\geq 2)})$ are Minkowski. Thus, by \cref{out ses}  $\Out(\Lambda)$ is Minkowski, as we wanted to see.
\end{proof}

The rest of this section is devoted to prove \cref{hyperbolic product asymmetric}. Recall that if $H$ is a torsion-free hyperbolic group, then $C_H(h)$ is cyclic for all non-trivial $h\in H$. In addition, if $H$ contains a normal abelian subgroup then $H$ is cyclic (see \cite[Part III $\Gamma$.3]{bridson2013metric}). Finally, if $M$ is a closed connected aspherical manifold of dimension $n\geq 3$ and $\pi_1(M)$ is hyperbolic then $\Out(\pi_1(M))$ is finite (see \cref{remark Out Minkowski}(2)). 

The arguments used in \cref{out hyperbolic finite} can be used to prove the next statement.

\begin{prop0}\label{product hyperbolic groups}
Let $M=M_1\times \cdots\times M_m$, where $M_i$ are closed aspherical manifolds such that $\pi_1(M_i)$ is hyperbolic and $\dim(M_i)\geq 3$. Then $\Out(\pi_1(M))$ is finite.
\end{prop0}

\begin{proof}
	Firstly, note that $\pi_1(M)=\pi_1(M_1)\times \cdots \times \pi_1(M_m)$. Let $H=\{f\in \Aut(\pi_1(M)):f(\pi_1(M_i))=\pi_1(M_i)\text{ for all $i$}\}$. As in \cref{out hyperbolic finite}, we will show that $H$ has finite index in  $\Aut(\pi_1(M))$. Since $\Inn(\pi_1(M))\trianglelefteq H$, this will imply that $[\Out(\pi_1(M)):H/\Inn(\pi_1(M))]<\infty$. 
	
	Let $e_i$ denote the trivial element of $\pi_1(M_i)$. We know that for every $\lambda=(\lambda_1,\dots,\lambda_m)\in\pi_1(M)$ we have $C_{\pi_1(M)}(\lambda)=\prod_{i=1}^{m}C_{\pi_1(M_i)}(\lambda_i)$. In addition, $C_{\pi_1(M_i)}(\lambda_i)=\pi_1(M_i)$ if $\lambda_i=e_i$ and $C_{\pi_1(M_i)}(\lambda_i)=\Z$ otherwise. As before, let $r(\lambda)$ denote the number of non-trivial entries of $\lambda$. We claim that that if $r(\lambda)=1$ and $f\in \Aut(\pi_1(M))$ then $r(f(\lambda))=1$. 
	
	Assume on the contrary, that $\lambda=(e_1,\dots,\lambda_i,\dots,e_m)$ and that $r({f(\lambda)})>1$. Since $f(C_{\pi_1(M)}(\lambda))=C_{\pi_1(M)}(f(\lambda))$ we can take the inverse morphism $f^{-1}:C_{\pi_1(M)}(f(\lambda))\longrightarrow C_{\pi_1(M)}(\lambda)$ and restrict it to $\Z^{r(f(\lambda))}\leq C_{\pi_1(M)}(f(\lambda))$. The morphism $\pi_i\circ f^{-1}_{|\Z^{r(f(\lambda))}}:\Z^{r(f(\lambda))}\longrightarrow \Z$ cannot be injective. If $a$ is a non trivial element of $\ker \pi_i\circ f^{-1}_{|\Z^{r(f(\lambda))}}$, then there exists a $j\neq i$ such that $\pi_j\circ f^{-1}_{|\langle a\rangle}:\langle a\rangle\longrightarrow \Lambda_j$ is injective. Since $\langle a\rangle\trianglelefteq C_{\pi_1(M)}(f(\lambda))$ and $\pi_j:C_{\pi_1(M)}(\lambda)\longrightarrow \pi_1(M_j)$ is surjective we can conclude that $\Z\cong\pi_j\circ f^{-1}_{|\langle a\rangle}(\langle a\rangle)\trianglelefteq \pi_1(M_j)$. But $\pi_1(M_j)$ is hyperbolic, so from the fact that it contains an abelian normal subgroup we can conclude $\pi_1(M_j)\cong \Z$. This is a contradiction with the fact that $\pi_1(M_j)$ is the fundamental group of a closed aspherical manifold of dimension $\dim(M_j)\geq 3$.
	
	In conclusion, any $f\in \Aut(\pi_1(M))$ permutes the factors of $\pi_1(M)$ and thus we can construct a group morphism to the permutation group of $m$ letters, $\phi:\Aut(\pi_1(M))\longrightarrow S_m$ such that $H=\ker \phi$. Consequently, $[\Aut(\pi_1(M)):H]<m!$.
	
	We are now ready to prove that $\Out(\pi_1(M))$ is finite. We proceed by induction on the number of factors $m$. If $m=1$ then $\Out(\pi_1(M))$ is finite (\cref{remark Out Minkowski}(2)). Assume now that $\Out(\pi_1(M_1)\times\cdots\times\pi_1(M_{m-1}))$ is finite. By \cref{out ses}, $\Out(\pi_1(M),\pi_1(M_m))$ is finite if $\Out(\pi_1(M_m))$, $\Out(\pi_1(M_1)\times\cdots\times\pi_1(M_{m-1}))$ and $H^1(\pi_1(M_1)\times\cdots\times\pi_1(M_{m-1}),\Zc\pi_1(M_m))$ are finite. But they are finite by induction hypothesis and the fact that $\Zc\pi_1(M_m)$ is trivial. Since $H\leq \Aut(\pi_1(M),\pi_1(M_m))\leq \Aut(\pi_1(M))$ and $[\Aut(\pi_1(M)):H]<\infty$ we have that $[\Out(\pi_1(M)):\Out(\pi_1(M),\pi_1(M_m))]<\infty$, which implies that $\Out(\pi_1(M))$ is finite.
\end{proof}

Note that the hyperbolic groups are centreless. Consequently, $\Zc\pi_1(M)$ is trivial, which implies that $M$ is almost asymmetric by \cref{main theorem1} and \cref{product hyperbolic groups}. In addition:

\begin{corollary0}
The group $\Aut(\pi_1(M))$ is Minkowski.
\end{corollary0}

\begin{proof}
Note that $\Zc\pi_1(M)$ is trivial, therefore there is a short exact sequence $$1\longrightarrow \pi_1(M) \longrightarrow \Aut(\pi_1(M))\longrightarrow \Out(\pi_1(M))\longrightarrow 1.$$
Since $ \Out(\pi_1(M))$ is finite and $\pi_1(M)$ is torsion-free, we conclude that $\Aut(\pi_1(M))$ is virtually torsion-free and hence Minkowski.
\end{proof}

\section{Combining the two cases: Aspherical locally homogeneous spaces}\label{sec:joining cases}

The aim of this section is to finish the proof of \cref{main theorem2} and to prove \cref{submain theorem2 discsym}. For the first task, the strategy is to combine the results on solvable Lie groups and semisimple Lie groups obtained in the previous sections in a similar way we proved \cref{out Minkowski with center}. 

Given a connected Lie group $G$, we would like to use the Levi decomposition $G=R\rtimes S$, where $S$ is the semisimple part and $R$ is the solvable radical, the maximal normal connected solvable Lie subgroup of $G$. However, $\Gamma\cap R$ is not a lattice in $R$ in general (see \cite{geng2015radicals}). 

Instead of using the solvable radical we will use the amenable radical, the maximal normal connected amenable Lie subgroup of $G$. Let us recall some basic properties of amenable groups.

\begin{lemma0}\cite[Proposition 4.1.6]{zimmer2013ergodic}\label{amenable group extensions}
Let $1\longrightarrow G_1 \longrightarrow G_2\longrightarrow G_3\longrightarrow 1$ be a short exact sequence of groups. Then $G_2$ is amenable if and only if $G_1$ and $G_3$ are amenable. 

Solvable Lie groups and compact Lie groups are amenable. In particular, finite groups and $\Z$ are amenable. 
\end{lemma0}

The lemma implies that virtually polycyclic groups are amenable. 

We can decompose $G=A\rtimes S_{nc}$ where $S_{nc}$ is semisimple with no compact factors and $A$ is the amenable radical. Moreover, $A=R\rtimes S_c$ where $S_c$ is compact and semisimple. If $\Gamma$ is a lattice in $G$ then $A\cap \Gamma$ is a lattice in $A$ (see \cite{gelander2020minimal,geng2015radicals}). In addition, a lattice $\Gamma$ of $G$ is amenable if and only if $G$ is amenable (see \cite[Proposition 4.1.11]{zimmer2013ergodic}). We have the short exact sequence
$$1\longrightarrow \Gamma\cap A\longrightarrow \Gamma \longrightarrow \Gamma/\Gamma\cap A\longrightarrow 1,$$
where the group $\Gamma\cap A$ is a lattice in the amenable radical $A$ and $\Gamma/\Gamma\cap A$ is a lattice in the semisimple Lie group without compact factors $S_{nc}$. We write $\Gamma\cap A=\Gamma_A$ and $\Gamma/\Gamma\cap A=\Gamma_{nc}$.

Our first goal is to see that $\Gamma_A$ is virtually polycyclic. \Cref{amenable lattice polycyclic by finite} is probably well-known to experts, but it is difficult to find a proof in the literature. Hence, we provide a proof for the sake of completeness.

\begin{lemma0}\label{amenable lattice polycyclic by finite}
Let $\Gamma$ be a lattice in an amenable group $A$. Then $\Gamma$ is virtually polycyclic.
\end{lemma0}

\begin{proof}
Denote by $\pi_{A/R}:A\longrightarrow A/R$ the quotient map and define $L=\overline{\pi_{A/R}(\Gamma)}^o$. It is connected solvable Lie group by \cite[8.24]{raghunathan2012discrete}. Moreover, $L$ is abelian since it is a connected solvable Lie subgroup of the compact Lie group $A/R$. Then, $\tilde{R}=\pi_{A/R}^{-1}(L)$ is a connected solvable group since $\tilde{R}/R$ is abelian and $R$ is connected.

$\tilde{R}\cap \Gamma$ is a lattice in $\tilde{R}$ (see \cite[Claim 2.2]{gelander2020minimal}). We claim that $\tilde{R}\cap \Gamma$ is polycyclic. Indeed, we have that $\tilde{R}/\tilde{R}\cap \Gamma$ is a compact solvmanifold. Then $\tilde{R}/\tilde{R}\cap \Gamma\cong R'/\Gamma'$, where $R'$ is the universal cover of $\tilde{R}$ and $\Gamma'$ is a lattice in a simply connected solvable Lie group thus it is polycyclic. Moreover, using the long exact sequence of homotopy groups for a fibration we obtain the short exact sequence $$1\longrightarrow \pi_1(\tilde{R})\cong \Z^n \longrightarrow \Gamma'\longrightarrow \tilde{R}\cap\Gamma \longrightarrow 1.$$
Since $\tilde{R}\cap\Gamma\cong \Gamma'/\Z^n$ we can conclude that $\tilde{R}\cap\Gamma$ is polycyclic. 

We want to show that $[\Gamma:\tilde{R}\cap \Gamma]<\infty$. Let $H=N_{A}(\tilde{R})$ and let $H^o$ be its connected component (then $R\leq\tilde{R}\leq H^o\leq H\leq A$). Thus $\Gamma\leq H$ and $|H/H^o|$ is bounded by \cite[Corollary $2.6$]{gelander2020minimal}. In addition, $H^o/\tilde{R}$ is a compact Lie group. 

In consequence, $[\Gamma:\Gamma\cap\tilde{R}]=[\Gamma:\Gamma\cap H^o][\Gamma\cap H^o:\Gamma\cap\tilde{R}]$. The first term is finite because $\Gamma\leq H$ and $|H/H^o|$ is bounded. The second term is finite because $\Gamma/\Gamma\cap\tilde{R}\cong \Gamma\tilde{R}/\tilde{R}=\pi_{H^o/\tilde{R}}(\Gamma)$ is a discrete subgroup in a compact Lie group and thus it is finite.
\end{proof}

Now that we know the structure of $\Gamma_A$ we are interested in the relation between $\Gamma_A$ and $\Gamma$. Since we want to study the automorphisms and outer automorphisms of $\Gamma$, it is natural to study whether $\Gamma_A$ is a characteristic subgroup of $\Gamma$. 

If $\Zc(S_{nc})\neq\{e\}$ then we consider the subgroup $\pi^{-1}(\Zc(\Gamma_{nc}))=\Gamma_A'$, which is a virtually polycyclic group since it fits in the short exact sequence $1\longrightarrow \Gamma_A\longrightarrow \Gamma_A'\longrightarrow \Zc(\Gamma_{nc})\longrightarrow 1$. In consequence, there is a short exact sequence $1\longrightarrow \Gamma_A'\longrightarrow \Gamma\longrightarrow \overline{\Gamma}_{nc}\longrightarrow 1 $, where $\overline{\Gamma}_{nc}$ is a lattice of the centreless connected semisimple Lie group $S_{nc}/\Zc(S_{nc})$. By \cref{facts semisimple lattices}(1), there exists a normal finite index torsion-free lattice $\Gamma_{nc}'\leq\overline{\Gamma}_{nc}$.  Then we can consider the short exact sequence $1\longrightarrow \Gamma_A'\longrightarrow \Gamma'\longrightarrow \Gamma_{nc}'\longrightarrow 1$. We will later prove that $\Out(\Gamma')$ is Minkowski. This fact, together with the previous short exact sequence, will imply that $\Out(\Gamma)$ is Minkowski by \cref{outer automorphism and finite index subgroups}.
\begin{lemma0}
$\Gamma_A'$ is a characteristic subgroup of $\Gamma'$. 
\end{lemma0} 

\begin{proof}

Assume that $\Gamma_A'$ is not characteristic. Thus, there exists $f\in \Aut(\Gamma')$ such that $f(\Gamma_A')$ has a non trivial projection $\Lambda$ in $\Gamma_{nc}'$. The group $\Lambda$ is a torsion-free virtually polycyclic subgroup and hence it has a characteristic subgroup $\Lambda'$ which is polycyclic (see \cite[\S 9.5]{lee2010seifert}). In particular, $\Lambda'$ is solvable and normal in the semisimple centreless lattice $\Gamma'_{nc}$, contradicting \cref{normal solvable groups lattices}. Consequently, $\Gamma_A'$ is characteristic.
\end{proof}

In consequence, $\Out(\Gamma')=\Out(\Gamma',\Gamma_A')$. \Cref{out ses} shows that to prove that $\Aut(\Gamma')$ and $\Out(\Gamma')$ are Minkowski it is enough to prove that $\Aut(\Gamma_A')$, $\Aut(\Gamma_{nc}')$, $Z^1_{\psi}(\Gamma_{nc}',\Zc(\Gamma_A'))$, $\Out(\Gamma_A')$, $\Out(\Gamma_{nc}')$ and $\overline{H}^1_\psi(\Gamma_{nc}',\Zc(\Gamma_A'))$ are Minkowski. We already know that $\Aut(\Gamma_A')$ and $\Out(\Gamma_A')$ are Minkowski because $\Gamma_A'$ is virtually polycyclic (\cref{outer automorphisms polycyclic}), and $\Out(\Gamma_{nc}')$ is Minkowski by \cref{out Minkowski with center}. Thus we only need to check that $\Aut(\Gamma_{nc}')$, $Z^1_{\psi}(\Gamma_{nc}',\Zc(\Gamma_A'))$, $\Out(\Gamma'_A)$ and $\overline{H}^1_\psi(\Gamma_{nc}',\Zc(\Gamma_A'))$ are Minkowski.

\begin{lemma0}
The group $\Aut(\Gamma_{nc}')$ is Minkowski. 
\end{lemma0}

\begin{proof}
Since $\Gamma_{nc}'$ is centreless, we have a short exact sequence
$$1\longrightarrow \Gamma_{nc}'\longrightarrow\Aut(\Gamma_{nc}')\longrightarrow \Out(\Gamma_{nc}')\longrightarrow 1.$$

Since $\Gamma_{nc}'$ is virtually torsion-free and $\Out(\Gamma_{nc}') $ is Minkowski, then $\Aut(\Gamma_{nc}')$ is Minkowski.
\end{proof}

\begin{lemma0}
The group $\overline{H}^1_\psi(\Gamma_{nc}',\Zc\Gamma_A')$ is Minkowski.
\end{lemma0}
\begin{proof}
Since $H^1_\psi(\Gamma_{nc}',\Zc\Gamma_A')$ is a finitely generated abelian group and there is a surjective morphism $H^1_\psi(\Gamma_{nc}',\Zc\Gamma_A')\longrightarrow \overline{H}^1_\psi(\Gamma_{nc}',\Zc\Gamma_A')$, the group $\overline{H}^1_\psi(\Gamma_{nc}',\Zc\Gamma_A')$ is finitely generated and abelian, hence it is Minkowski.
\end{proof}

\begin{lemma0}
The group $Z^1_{\psi}(\Gamma'_{nc},\Zc(\Gamma'_A))$ is Minkowski.
\end{lemma0}

\begin{proof}
There is a short exact sequence $$1\longrightarrow \overline{B}^1_{\psi}(\Gamma'_{nc},\Zc(\Gamma'_A))\longrightarrow Z^1_{\psi}(\Gamma'_{nc},\Zc(\Gamma'_A))\longrightarrow\overline{H}^1_\psi(\Gamma_{nc}',\Zc\Gamma_A')\longrightarrow 1.$$ 

We have that $\overline{B}^1_{\psi}(\Gamma'_{nc},\Zc(\Gamma'_A))\cong \Zc(\Gamma'_A)/\Zc(\Gamma')$ (see \cref{out ses}), therefore $\overline{B}^1_{\psi}(\Gamma'_{nc},\Zc(\Gamma'_A)) $ is a finitely generated abelian group and hence it is Minkowski. Since $\overline{H}^1_\psi(\Gamma_{nc}',\Zc\Gamma_A')$ is also a finitely generated abelian group, we use \cref{Minkowski ses} to conclude that $Z^1_{\psi}(\Gamma'_{nc},\Zc(\Gamma'_A))$ is Minkowski.
\end{proof}

These three lemmas complete the proof of \cref{main theorem2}.

\begin{remark0}
Note that we can use the proof of \cref{main theorem2} to conclude that if $\Gamma$ is a torsion-free group which fits in a short exact sequence $1\longrightarrow \Gamma_A\longrightarrow \Gamma\longrightarrow \Gamma_{nc}\longrightarrow1$ where $\Gamma_A$ is virtually polycyclic and $\Gamma_{nc}$ is a lattice in a centreless semisimple Lie group then $\Out(\Gamma)$ is Minkowski. Not all groups of this form are lattices in connected Lie groups, see \cite[Theorem 7.5]{BauesKamishima}. Baues and Kamishima introduced in \cite{BauesKamishima} the notion of the closed aspherical manifolds with large symmetry. From the definition of these manifolds \cite[Definition 1.1]{BauesKamishima} one can deduce that their fundamental group fulfils the hypothesis of this remark and hence the outer automorphism group of their fundamental group is Minkowski.  
\end{remark0}

\begin{remark0}
The homeomorphism group of a coset space $G/\Gamma$ where $G$ is a connected Lie group and $\Gamma$ is a cocompact lattice is not necessarily Jordan. Indeed, we can take a non-compact semisimple Lie group $G'$ and a cocompact lattice $\Gamma'$ such that a connected compact maximal subgroup $K'$ containing $\SU(2)$ acts effectively on $G'/\Gamma'$. By \cite{riera2017non} we know that  $\Homeo(G'/\Gamma'\times T^2)$ is not Jordan. Hence, taking $G=G'\times \R^2$ and $\Gamma=\Gamma'\times \Z^2$ we obtain a coset space $G/\Gamma$ whose homeomorphism group is not Jordan. On the other hand, A. Golota proved in \cite{golota2023finite} that if $G$ is a complex Lie group and $\Gamma$ is a cocompact lattice, then the group of automorphisms $\Aut(G/\Gamma)$ is Jordan.
\end{remark0}

The end of this section is devoted to prove \cref{submain theorem2 discsym}. We will need the following two results from \cite{lee2010seifert}.

\begin{thm0}\label{tor-sym homotopy equivalent}\cite[Theorem 11.7.29]{lee2010seifert}
Let $M$ be a closed connected aspherical manifold whose fundamental group fits into a short exact sequence $1\longrightarrow \Gamma_A\longrightarrow \pi_1(M)\longrightarrow \Gamma_{nc}\longrightarrow 1$, where $\Gamma_A$ is virtually polycyclic and $\Gamma_{nc}$ is a cocomapct lattice in a semisimple connected Lie group without compact factors. Then there exists a closed connected aspherical manifold $M'$ which is homotopically equivalent to $M$ and satisfies $\A(M')=\rank\Zc\pi_1(M)$.
\end{thm0}

\begin{thm0}\label{tor-sym solvable radical R}\cite[Theorem 11.7.28]{lee2010seifert}
Let $H\setminus G/\Gamma$ be a closed aspherical locally homogeneous space with $G$ simply connected and let $R$ be the solvable radical of $G$. Assume that $R\cap \Gamma$ is a lattice in $R$ such that any automorphism of $R\cap \Gamma$ can be extended to $R$, and that $\exp:\mathcal{L}(R)\longrightarrow R$ is surjective. Then $\A(H\setminus G/\Gamma)=\rank\Zc\Gamma$.
\end{thm0}

\begin{proof}[Proof of \cref{submain theorem2 discsym}]
Note that $\Out(\Gamma)$ is Minkowski by \cref{main theorem2}, hence $\D(H\setminus G/\Gamma)\leq\rank\Zc\Gamma$ by \cref{main theorem1}. On the other hand, it is clear that $\A(H\setminus G/\Gamma)\leq \D(H\setminus G/\Gamma)$. Therefore, it suffices to show that $\A(H\setminus G/\Gamma)=\rank\Zc\Gamma$ to conclude that $\D(H\setminus G/\Gamma)=\rank\Zc\Gamma$. 

By \cref{amenable lattice polycyclic by finite} $\Gamma$ satisfies the conditions of \cref{tor-sym homotopy equivalent}, hence there exists a closed connected aspherical manifold $M'$ homotopically equivalent to $H\setminus G/\Gamma$ such that $\A(M')=\rank\Zc\Gamma$. On the other hand, the Farrel-Jones conjecture is true for cocompact lattices in connected Lie groups (see \cite{bartels2012borel,kammeyer2016farrell}), which implies the Borel conjecture for closed aspherical locally homogeneous space of dimension equal or greater than 5 (see \cite[Proposition 0.3 (ii)]{bartels2012borel}). In consequence, if $\dim(H\setminus G/\Gamma)\geq 5$ then $H\setminus G/\Gamma\cong M'$ and hence $\A(H\setminus G/\Gamma)=\rank\Zc\Gamma$. Thus, it remains to study the cases where $\dim(H\setminus G/\Gamma)\leq 4$.

If $\dim(H\setminus G/\Gamma)=1$ or $2$ then $\A(H\setminus G/\Gamma)=\rank\Zc\Gamma$ by classical results and if $\dim(H\setminus G/\Gamma)=3$ then $\A(H\setminus G/\Gamma)=\rank\Zc\Gamma$ by \cite[Corollary 8.3]{gabai1992convergence} and \cite[Theorem 1.1]{casson1994convergence}. Thus, it only remains to study the case where $\dim(H\setminus G/\Gamma)=4$. 

First, assume that the noncompact semisimple part $S_{nc}$ of $G$ is trivial. Then $\Gamma$ is virtually polycyclic. Virtually polycyclic groups are Freedman good (this property was called "good" and it was introduced in \cite{freedman1983disk}). If the fundamental group is Freedman good, the Farrel-Jones conjecture implies the Borel conjecture in dimension $4$ (see \cite[Proposition 0.3 (ii)]{bartels2012borel}) and hence $\A(H\setminus G/\Gamma)=\rank\Zc\Gamma$. 

In consequence, it only remains to study the case where $S_{nc}\neq \{e\}$. We can assume two extra hypothesis without losing generality. Firstly, let $q:\tilde{G}\longrightarrow G$ be the universal cover of $G$. Let $\tilde{H}$ be a maximal compact connected subgroup of $\tilde{G}$ inside $q^{-1}(H)$. Then $q^{-1}(\Gamma)=\tilde{\Gamma}$ is a cocompact lattice of $\tilde{G}$ and $H\setminus G/\Gamma\cong \tilde{H}\setminus\tilde{G}/\tilde{\Gamma}$. Thus, we can assume that $G$ is simply connected and hence its solvable radical $R$ and its semisimple part $S$ are also simply connected. Secondly, note that if $K$ is a compact connected normal subgroup of $G$ and $p:G\longrightarrow G/K$ is the quotient map then $\Gamma\cap K=\{e\}$ since $\Gamma$ is torsion-free. This implies that $H\setminus G/\Gamma\cong p(H)\setminus p(G)/\Gamma$. Therefore, we can assume that $G$ has no compact factors.  

Note that $\dim(R)+\dim(S)-\dim(H)=4$ and since $S_{nc}\neq\{e\}$ then $\dim(S)-\dim(H)\geq 2$. Hence, $\dim(R)\leq 2$. There are four simply connected solvable groups of dimension less or equal than $2$: $\dim(R)=0$ and $R=\{e\}$, $\dim(R)=1$ and $R\cong \R$, and $\dim(R)=2$ and $R\cong \R^2$ or $R\cong \Aff(\R)^0\cong\R\rtimes\R$. To construct the desired torus action on $H\setminus G/\Gamma$ with $R$ belonging to one of this four cases we will use \cref{tor-sym solvable radical R}. We need to check that the hypothesis of \cref{tor-sym solvable radical R} are satisfied in these four cases. Note that in all the cases the exponential map is surjective and any lattice automorphism extends uniquely to an automorphism of the group. This is because $R$ is either abelian or solvable of type (R) (see \cite[\S 6.3]{lee2010seifert}). It only remains to see that $\Gamma\cap R$ is a lattice in $R$. 

Assume that $R$ is abelian, then if $C$ is a compact factor of $S$ acting trivially on $R$ then $C$ is normal on $G$ and therefore $C=\{e\}$. Therefore no compact factor of $S$ acts trivially on $R$. This implies by \cite[Theorem 1.3 (i)]{geng2015radicals} that $R\cap\Gamma$ is a lattice in $R$. Thus, $\A(H\setminus G/\Gamma)=\rank\Zc\Gamma$ by \cref{tor-sym solvable radical R}.

If $R\cong \Aff(\R)^0$ then $\Aut(R)\cong R$ and hence any compact factor $C$ of $S$ needs to act trivially on $R$. Consequently, $C\trianglelefteq G$ and therefore $C=\{e\}$. By \cite[Theorem 1.3 (i)]{geng2015radicals} if $\Gamma$ were a lattice of $G$ then $R\cap\Gamma$ would be a lattice of $R$. But $R$ does not admit any lattice (see \cite[pg. 82]{bock2016low}). Thus, a closed aspherical locally homogeneous 4-manifold with solvable radical $R$ does not exists. 
\end{proof}

\section{Jordan property and cohomology}

This section is devoted to proving \cref{aspherical Jordan cohomology}. We start with the following general lemma.

\begin{lemma0}\label{submain theorem2}
	Let $\Gamma$ be a finitely generated group and let $\Zc_i\Gamma$ be the $i$-th term of the upper central series. Assume that $\Gamma/\Zc_i\Gamma$ is finitely generated, centreless and that $\Out(\Gamma/\Zc_i\Gamma)$ is Minkowski for some $i$. Then $\Out(\Gamma)$ is Minkowski.
\end{lemma0}

\begin{proof}
	Recall that the upper central series $\{e\}=\Zc_0\Gamma\trianglelefteq \Zc_1\Gamma\trianglelefteq \Zc_2\Gamma\trianglelefteq\cdots$ of a group $\Gamma$ is a subnormal series where $\Zc_{i}\Gamma/\Zc_{i-1}\Gamma\cong \Zc(\Gamma/\Zc_{i-1})$ for all $i\geq 0$. The group $\Gamma$ is nilpotent if and only if $\Gamma\cong \Zc_i\Gamma$ for some $i$. Moreover, the groups $\Zc_i\Gamma$ are characteristic and each group morphism $\psi_i:\Gamma/\Zc_i\Gamma\longrightarrow \Out(\Zc_i\Gamma)$ is trivial. 
	
	Consequently, the group $\Out(\Gamma)$ is Minkowski if $\Out(\Zc_i\Gamma)$, $\Out(\Gamma/\Zc_i\Gamma)$ and $H^1(\Gamma/\Zc_i\Gamma,\Zc\Gamma_i/\Zc\Gamma_{i-1})$ are Minkowski. But they are Minkowski by hypothesis and \cref{outer automorphisms polycyclic}, obtaining the desired conclusion. 
\end{proof}

If we assume that $\Gamma$ is torsion-free and $\Gamma/\Zc_i\Gamma$ acts properly on a contractible manifold $\tilde{X}$ then we can use the Seifert fiber construction in \cite[Chapter 7]{lee2010seifert} to construct closed connected aspherical manifolds $M$ such that $\pi_1(M)\cong \Gamma$. In addition, if $\tilde{X}/(\Gamma/\Zc_i\Gamma)$ is a closed aspherical manifold, then $M$ can be seen as an iterated principal torus bundle over $\tilde{X}/(\Gamma/\Zc_i\Gamma)$. A special case of \cref{submain theorem2} and \cref{main theorem1} is the next corollary. 

\begin{corollary0}\label{product manifolds minkowski}
	If $M$ is a closed connected aspherical manifold such that $\Zc\pi_1(M)$ is trivial and $\Out(\pi_1(M))$ is Minkowski, then $\Homeo(M\times T^n)$ is Jordan.
\end{corollary0}

Finally, there exists closed hyperbolic $3$-manifolds $N$ which are integral homology spheres (see \cite{benedetti1992lectures,thurston2022geometry}). We have that $\Zc(\pi_1(N))$ is trivial and $\Out(\pi_1(N))$ is finite by \cite[5.4 A]{gromov1987hyperbolic}. Thus, $H^*(N\times T^2)\cong H^*(S^3\times T^2)$, $\Homeo(N\times T^2)$ is Jordan by \cref{product manifolds minkowski}  and $\Homeo(S^3\times T^2)$ is not Jordan by \cite{riera2017non}. This completes the proof of \cref{aspherical Jordan cohomology}. 

Another example can be produced using Brieskorn manifolds $\Sigma(p,q,r)$. Recall that $\Sigma(p,q,r)=\{(z_1,z_2,z_3)\in\C^3:z_1^p+z_2^q+z_3^r=0, \sum_{i=1}^3|z_i|^2=1\}$. If $p$, $q$ and $r$ are relatively prime and $\tfrac{1}{p}+\tfrac{1}{q}+\tfrac{1}{r}<1$ then $\Sigma(p,q,r)$ is a closed connected aspherical integral homology $3$-sphere (see \cite{milnor19753}). Moreover, $\Sigma(p,q,r)$ are Seifert manifolds with $\Zc\pi_1(\Sigma(p,q,r))\cong \Z$ (see \cite[\S 14.11]{lee2010seifert}), hence there exists a $S^1$-action on $\Sigma(p,q,r)$ inducing the short exact sequence
$$1\longrightarrow \Z\longrightarrow \pi_1(\Sigma(p,q,r))\longrightarrow \Inn\pi_1(\Sigma(p,q,r))\longrightarrow 1$$
where $\Inn\pi_1(\Sigma(p,q,r))$ is centreless. In addition, $\Inn\pi_1(\Sigma(p,q,r))$ is a subgroup of isometries of the hyperbolic plane. Hence, it contains a centreless torsion-free Fuchsian subgroup $Q$ of finite index. Since $\Out(Q)$ is virtually torsion-free then $\Out(\Inn\pi_1(\Sigma(p,q,r)))$ is also virtually torsion-free (see \cite[Lemma 2.4, Corollary 2.6]{metaftsis2006residual}) and therefore $\Out(\Inn\pi_1(\Sigma(p,q,r)))$ is Minkowski.

The manifold $\Sigma(p,q,r)\times T^2$ is a closed connected aspherical manifold such that $H^*(\Sigma(p,q,r)\times T^2)\cong H^*(S^3\times T^2)$. Moreover, since $\Zc(\pi_1(\Sigma(p,q,r)\times T^2))\cong\Z^3$ we have a central extension 
$$1\longrightarrow \Z^3\longrightarrow \pi_1(T^2\times \Sigma(p,q,r))\longrightarrow \Inn\pi_1(\Sigma(p,q,r))\longrightarrow 1.$$

By \cref{submain theorem2}, the group $\Out(\pi_1(T^2\times \Sigma(p,q,r)))$ is Minkowski. Therefore, \cref{main theorem1} implies that $\Homeo(T^2\times \Sigma(p,q,r)))$ is Jordan.

\begin{remark0} A lot of the results about the Jordan property on $\Homeo(M)$ of a closed connected manifold $M$ rely on the cohomology of $M$, for example in \cite{mundet2019finite} its is proven that $\Diff(M)$ is Jordan if $M$ is an integral homology sphere or $\chi(M)\neq 0$. However, \cref{aspherical Jordan cohomology} shows that the Jordan property on $\Homeo(M)$ and $\Diff(M)$ does not only depend on the cohomology of $M$ in general. 
\end{remark0}

\begin{remark0}
	We also note that $\D(N\times T^2)<\D(\Sigma(p,q,r)\times T^2)<\D(S^3\times T^2)$, thus all three manifolds have different discrete degree of symmetry.
	
	We have $\D(N\times T^2)=2$, since $\D(N\times T^2)\leq 2$ by \cref{main theorem1} and $N\times T^2$ admits a $T^2$ free action. Similarly, $\D(\Sigma(p,q,r)\times T^2)=3$, since $\D(\Sigma(p,q,r)\times T^2)\leq 3$ by \cref{main theorem1} and $\Sigma(p,q,r)\times T^2$ admits an action of $T^3$. Finally, $T^4$ acts effectively on $S^3\times T^2$, hence $\D(S^3\times T^2)\geq 4$. 
\end{remark0}


\section{Remarks when the discrete degree of symmetry is close to the dimension of the aspherical manifold: Proof of \cref{submain theorem 1}}\label{sec discsym}

We have seen that if $M$ is a closed $n$-dimensional aspherical manifold with $\Zc\pi_1(M)$ finitely generated and $\Out(\pi_1(M))$ Minkowski then $\D(M)=n$ if and only if $M\cong T^n$. An interesting question is whether there exist similar rigidity results when $\D(M)$ is close to $n$, for example $\D(M)=n-1$. The aim of this section is to answer this question provided we have an additional hypothesis on $\pi_1(M)$.

We start by recalling what happens in low dimensions. Let $M$ be a $2$-dimensional closed connected aspherical manifold. If $M$ is orientable then $M$ is either a torus $T^2$ and $\A(M)=2$ or $M$ is a surface $\Sigma_g$ of genus $g\geq 2$ and $\A(M)=0$. If $M$ is not orientable, then $M$ is either the Klein bottle $K$ and $\A(M)=1$ or $M$ has a surface $\Sigma_g$ of genus $g\geq 2$ as an orientable 2-cover, hence $\A(M)=0$.

If $M$ is a closed connected $3$-dimensional aspherical manifold with an effective $S^1$ action, then $M$ falls in one of the following $4$ cases (see \cite[\S 14.4]{lee2010seifert}):
\begin{itemize}
	\item[1.] $M\cong T^3$.
	\item[2.] $M$ is homeomorphic to $K\times S^1$ or $SK$, the non-trivial principal $S^1$-bundle over $K$.
	\item[3.] $M\cong H/\Gamma$, where $H$ is the $3$-dimensional Heisenberg group and $\Gamma $ is a lattice of $H$.
	\item[4.] $\Zc\pi_1(M)\cong \Z$ and $\Inn \pi_1(M)\cong \pi_1(M)/ \Zc\pi_1(M)$ is centreless.
\end{itemize}

Note that in all cases, we have a central extension $$1\longrightarrow \Z\longrightarrow \pi_1(M)\longrightarrow Q\longrightarrow 1$$ where $Q$ acts effectively, properly and cocompactly on $\R^2$.

In the first case $\A(M)=3$ and in the third and fourth cases $\A(M)=1$. In the second case, it is clear that $\A(K\times S^1)=2$, so we focus on $SK$. One can see that $\pi_1(SK)\cong \Z^2\rtimes_\phi \Z$, where $\phi:\Z\longrightarrow\Gl(2,\Z)$ satisfies that
\[
\phi(1)=\left(\begin{matrix}
	0 & 1 \\
	1 & 0
\end{matrix}\right).
\]

Since $\phi(1)^2=Id$, we have that $SK$ is a flat solvmanifold and therefore $\A(SK)=\rank\Zc\pi_1(SK)$.
\begin{lemma0}
	Let $\phi:\Z\longrightarrow \Gl(n,\Z)$ be group morphism such that $\phi(1)$ has finite order $a$. Then $\Zc(\Z^n\rtimes_\phi \Z)=\operatorname{Fix}(\phi)\times a\Z$, where $\operatorname{Fix}(\phi)=\{v\in \Z^n:\phi(1)v=v\}$.
\end{lemma0}

\begin{proof}
A straightforward computation shows that given $(v,t),(w,s)\in\Z^n\rtimes_\phi \Z$ we have $(v,t)(w,s)=(w,s)(v,t)$ if and only if $(Id-\phi(s))v=(Id-\phi(t))w$. If $(v,t)\in\Zc(\Z^n\rtimes_\phi \Z)$ then we have $(Id-\phi(t))w=0$ for all $w\in \Z^n$, by taking $s=0$. Thus, $\phi(t)=Id$ and hence $t\in a\Z$. 

If we assume that $w=0$, then $(Id-\phi(s))v=0$ for all $s\in \Z$. Thus, $v\in \operatorname{Fix}(\phi)$. Consequently, $\Zc(\Z^n\rtimes_\phi \Z)\subseteq\operatorname{Fix}(\phi)\times a\Z$. The other inclusion follows from the fact that if $(v,t)\in \operatorname{Fix}(\phi)\times a\Z$ then $(Id-\phi(s))v=0=(Id-\phi(t))w$ for any $s\in \Z$ and $w\in \Z^n$.
\end{proof}

Applying the previous lemma to $\pi_1(SK)$, we obtain $\operatorname{Fix}(\phi)=\langle(1,1)\rangle$ and therefore $\Zc\pi_1(SK)\cong \Z^2$ and $\A(SK)=2$. Thus, $K$ and $SK$ are the two only $3$-dimensional aspherical manifolds such that where $\A(M)=\dim(M)-1$. The next proposition generalizes the previous fact to arbitrary dimension. It is probably well-known to experts but we could not find a proof in the literature, so we provide it for the sake of completeness.

\begin{prop0}\label{absym=n-1}
	Let $M$ be a closed aspherical manifold of dimension $n$. If $\A(M)=n-1$ then $M\cong K\times T^{n-2}$ or $M\cong T^{n-3}\times SK$. In particular, $M$ is always a non-orientable flat solvmanifold.
\end{prop0}

\begin{proof}
	Let $H\leq T^{n-1}$ be the isotropy subgroup of the principal orbit of the action. Since  $T^{n-1}$ is abelian and its action on $M$ is effective we can conclude that $H$ is trivial. Therefore, the principal orbits of the action have dimension $n-1$. In this case we say that we have a cohomogenity one action. The next theorem describes the cohomogenity one actions.
	
	\begin{thm0}\cite[Theorem A]{galaz2018cohomogeneity}
		Let $M$ be a $n$-dimensional closed connected manifold with a cohomogenity one action of a compact Lie group $G$. Let $H$ be the isotropy subgroup of a principal orbit. Then we have one of these two options:
		\begin{itemize}
			\item[1.] The quotient $M/G\cong S^1$. Then $M$ is equivariantly homeomorphic to the total space of a fiber bundle $\pi: M\longrightarrow S^1$ is a fiber bundle with fiber $G/H$. The action does not have exceptional orbits.
			\item[2.] The quotient $M/G\cong [-1,1]$. Then $M$ is the union of
			two fiber bundles over the two singular orbits with isotropy subgroups $K_+$ and $K_-$ whose fibers are cones
			over spheres or the Poincaré homology sphere. More explicitly,
			$$M = G\times_{K_-}C(K_-/H)\cup_{G/H} G\times_{K_+}C(K_+/H)$$
			where $C(K_{\pm}/H)$ denotes the cone over $K_{\pm}/H$, which are spheres or Poincaré homology spheres.  The exceptional orbits $G/K_\pm$ correspond to the preimages of $\pm1\in [-1,1]$.
		\end{itemize}
	\end{thm0}
	If the cohomogenity one action of $T^{n-1}$ on $M$ does not have exceptional orbits, then the action is free and the orbit map $\pi:M\longrightarrow S^1$ is a principal $T^{n-1}$-bundle. Since the base is $S^1$ the principal bundle is trivial and therefore $M\cong T^n$ and $\A(M)=n$, which is not possible. Therefore the action has exceptional orbits. 
	
	Since $M$ is aspherical, the evaluation map $ev_x:T^{n-1}\longrightarrow M$ such that $ev_x(g)=gx$ for all $g\in T^{n-1}$ induces an injective group morphism ${ev_x}_*:\pi_1(T^{n-1})\longrightarrow \pi_1(M)$ for any $x\in M$ (\cite[Lemma 3.1.11]{lee2010seifert}). Consequently, all the isotropy subgroups of the action are discrete. The quotients $K_\pm/H\cong K_\pm$ are homeomorphic to a sphere or a Poincaré homology sphere. Since they are also discrete we obtain that $K_\pm\cong S^0\cong \Z/2$. Then 
	$$M = T^{n-1} \times_{\Z/2}I_+\cup_{T^{n-1}}  T^{n-1} \times_{\Z/2} I_-,$$
	where $I_+=[0,1]$ and $I_-=[-1,0]$. Let $Z_+=\pi_1(T^{n-1}\times_{\Z/2}I_+)$ and $Z_-=\pi_1(T^{n-1}\times_{\Z/2}I_-)$. The principal $\Z/2$-bundles $ T^{n-1} \times I_\pm\longrightarrow  T^{n-1} \times_{\Z/2}I_\pm$ induce two short exact sequences of fundamental groups $$ 1\xrightarrow{}\Z^{n-1}\xrightarrow{i_\pm}Z_\pm\xrightarrow{}\Z/2\xrightarrow{}1.$$
	Note that $Z_+$ and $Z_-$ are isomorphic to $\Z^{n-1}$. By the Seifert-van Kampen theorem, $\pi_1(M)=Z_+*_{\Z^{n-1}}Z_-$ where the amalgamated product is induced by the inclusions $i_\pm:\Z^{n-1}\longrightarrow Z_\pm$, which we are going to describe explicitly. 
	
	Firstly, there exist primitive elements $\alpha_\pm\in Z_\pm$ such that $\alpha_\pm\notin i_\pm(\Z^{n-1})$. Furthermore, there exist two $x_\pm\in \Z^{n-1}$ such that $i_\pm(x_\pm)=2\alpha_\pm$. We have two possibilities, that $x_+=x_-$ or that $x_+\neq x_-$.
	
	Firstly, assume that $x_+=x_-$. If we denote $x_+=x_-=x$ then we can choose generators $\{x,y_1,...,y_{n-2}\}$ of $\Z^{n-1}$ such that $\{\alpha_+,i_+(y_1),...,i_+(y_{n-2})\}$ generates $Z_+$ and $\{\alpha_-,i_-(y_1),...,i_-(y_{n-2})\}$ generates $\Z_-$. Thus,
	$$\pi_1(M)\cong \langle\alpha_+,\alpha_-|\alpha_+^2=\alpha_-^2\rangle\times \Z^{n-2}\cong \pi_1(K)\times \Z^{n-2}=\pi_1(K\times T^{n-2}).$$ 
	Since $\pi_1(M)$ is the fundamental group of a flat manifold and the Borel conjecture holds for these groups (see \cite{bartels2012borel}) we conclude that $M\cong K\times T^{n-2}$.
	
	We now assume that $x_+\neq x_-$. As before, we can choose a generator set $\{x_+,x_-,y_1,...,y_{n-3}\}$ of $\Z^{n-1}$ such that $\{\alpha_+,i_+(x_-),i_+(y_1),...,i_+(y_{n-3})\}$ generates $Z_+$ and $\{i_-(x_+),\alpha_-,i_-(y_1),...,i_-(y_{n-3})\}$ generates $Z_-$. Thus, $$\pi_1(M)\cong \langle\alpha_+,i_+(x_-),i_-(x_+),\alpha_-|\alpha_+^2= i_-(x_+), \alpha_-^2=i_+(x_-),[\alpha_+,i_+(x_-)]=1,[\alpha_-,i_-(x_+)]=1\rangle\times \Z^{n-3}.$$
	The presentation of the first factor $\Lambda$ can be rearranged to obtain a new presentation
	\begin{align*}
		\Lambda \cong & \langle a,b| [a,b^2]=1,[a^2,b]=1\rangle\\
		\cong & \langle a,b,c,d|c=ab,d=ba, ac=bd,cb=da\rangle\\
		\cong & \langle a,c,d|[c,d]=1,aca^{-1}=d,ada^{-1}=c\rangle.
	\end{align*}
	With the last presentation, we can define an isomorphism $f:\Lambda\longrightarrow \Z^2\rtimes_\phi\Z$ such that $f(a)=(0,0,1)$, $f(b)=(1,0,0)$ and $f(c)=(0,1,0)$. Consequently, $\pi_1(M)\cong\pi_1(SK\times T^{n-3})$. Since $\pi_1(M)$ is the fundamental group of a flat manifold and the Borel conjecture holds for these groups we conclude that $M\cong SK\times T^{n-3}$.
	
\end{proof}

\begin{corollary0}
	Let $M$ be a closed connected aspherical manifold of dimension $n$. If $M\ncong T^n$ and $M$ is orientable then $\A(M)\leq n-2$.
\end{corollary0}

\begin{remark0}
	The hypothesis of $M$ being aspherical is essential. For example $\A(S^2)=1$ and $S^2\ncong K\times T^{n-2},T^{n-3}\times SK$.
\end{remark0}

\begin{thm0}\label{discsym=n-1}
	Let $M$ be a closed connected aspherical manifold such that $\Zc(\pi_1(M))$ is finitely generated and $\Out(\pi_1(M))$ is Minkowski. Assume that $\Inn(\pi_1(M))$ has an element of infinite order. Then $\D(M)=n-1$ if and only if $M\cong K\times T^{n-2}$ or $M\cong T^{n-3}\times SK$.
\end{thm0}

\begin{proof}
	Because of the hypothesis on $\pi_1(M)$ we know that $\Zc(\pi_1(M))\cong \Z^{n-1}$. Let $x\in \Inn\pi_1(M)$ be an element of infinite order and consider the commutative diagram
	\[
	\begin{tikzcd}
		1\ar{r}{} & \Zc\pi_1(M)\ar{d}{Id}\ar{r}{} & p^{-1}(\langle x\rangle)\ar{d}{}\ar{r}{} & \langle x\rangle\ar{d}{}\ar{r}{} & 1\\
		1\ar{r}{} & \Zc\pi_1(M)\ar{r}{} & \pi_1(M)\ar{r}{p} & \Inn\pi_1(M)\ar{r}{} & 1
	\end{tikzcd}
	\]
	where the vertical arrows are inclusion maps. The upper short exact sequence is also central and therefore it is classified by $H^2(\Z,\Z^{n-1})=0$. In consequence, $p^{-1}(\langle x\rangle)\cong \Z^n$. We consider now the covering $q:\tilde{M}/p^{-1}(\langle x\rangle)\longrightarrow M$. Note that $H^*(\tilde{M}/p^{-1}(\langle x\rangle),\Z)\cong H^*(\Z^n,\Z)$ and therefore $H^n(\tilde{M}/p^{-1}(\langle x\rangle),\Z)$ is not trivial. Consequently, $\tilde{M}/p^{-1}(\langle x\rangle)$ is a closed connected aspherical manifold and $ \tilde{M}/p^{-1}(\langle x\rangle)\cong T^n$. Since $M$ and $\tilde{M}/p^{-1}(\langle x\rangle)$ are compact, we obtain that $q$ is a finite covering map and therefore $p^{-1}(\langle x\rangle)\leq \pi_1(M)$ has finite index and $\pi_1(M)$ is the fundamental group of a flat manifold. Thus, we have that $\D(M)=\A(M)$ by \cref{submain theorem2 discsym} and the result follows from \cref{absym=n-1}.
\end{proof}

It is an interesting question to know whether $\Inn\pi_1(M)$ always has an element of infinite order when $M$ is closed aspherical manifold not homeomorphic to a torus and such that $\Zc\pi_1(M)$ is finitely generated. Note that $\Inn\pi_1(M)$ cannot be finite by \cref{schur inner theorem} unless $M$ is homeomorphic to a torus. Thus, if $\Inn\pi_1(M)$ does not contain elements of infinite order then $\pi_1(M)$ is abelian or $\Inn\pi_1(M)$ is infinite periodic. 

We present some evidences that support that the answer to the question is affirmative. Firstly, let $\Gamma$ and $\Lambda$ be two groups. Assume that we have a surjective group morphism $p:\Gamma\longrightarrow \Lambda$. Then there is an induced surjective group morphism $p':\Inn\Gamma\longrightarrow \Inn\Lambda$ which sends a conjugation $c_\gamma$ to $c_{p(\gamma)}$. Thus, if $\Inn\Lambda$ has an element of infinite order, then $\Inn\Gamma$ also has an element of infinite order. 

Assume now that we have an inclusion $i:\Gamma\longrightarrow \Lambda$ instead. Then there is an induced injective group morphism $i':\Inn\Gamma\longrightarrow \Inn\Lambda$ which sends a conjugation $c_\gamma$ to $c_{i(\gamma)}$. In particular, if $\Inn\Gamma$ has an element of infinite order, then $\Inn\Gamma$ also has an element of infinite order. We can deduce the following corollary from these observations:

\begin{corollary0}
	Let $M$ be a closed aspherical manifold, then:
	\begin{itemize}
		\item[(1)] If $M'\longrightarrow M$ is a covering and $\Inn\pi_1(M')$ has an element of infinite order then $\Inn\pi_1(M)$ has an element of infinite order.
		\item[(2)] Suppose that we have a fibration of closed connected aspherical manifolds $M'\longrightarrow M\longrightarrow M''$. If $\Inn\pi_1(M')$ or $\Inn\pi_1(M'')$ have an element of infinite order, then $\Inn\pi_1(M)$ has an element of infinite order.
	\end{itemize}
\end{corollary0}

\begin{proof}
	In the first case, we use that $\pi_1(M')\leq \pi_1(M)$ and in the second case the short exact sequence $1\longrightarrow\pi_1(M')\longrightarrow \pi_1(M)\longrightarrow \pi_1(M'')\longrightarrow 1$ together with the observations above.
\end{proof}

Note that if $\Gamma$ is a non-abelian polycyclic group or $\Gamma$ is torsion-free centreless then $\Inn\Gamma$ has elements of infinite order. Thus, closed aspherical locally homogeneous spaces and closed aspherical manifolds whose fundamental group is hyperbolic have an infinite order element in the inner automorphism group of their fundamental group. Moreover, any closed aspherical manifold $M$ not homeomorphic to a torus constructed using fibrations of these two classes of aspherical manifolds will have an element of infinite order in $\Inn\pi_1(M)$.

There are cohomological restrictions to $\Inn\pi_1(M)$ being infinite periodic. Assume that $M$ is a closed connected aspherical manifold such that $\Inn\pi_1(M)$ is infinite periodic. Then $H^1(\Inn\pi_1(M),\Z)=\Hom(\Inn\pi_1(M),\Z)$ is trivial and the inflation-restriction exact sequence becomes

$$ 1\longrightarrow H^1(\pi_1(M),\Z)\longrightarrow H^1(\Zc\pi_1(M),\Z) \longrightarrow H^2(\Inn\pi_1(M),\Z)\longrightarrow H^2(\pi_1(M),\Z)$$

where we see $\Z$ as a trivial $\pi_1(M)$-module. In particular, if $\rank H^1(M,\Z)>\rank \Zc\pi_1(M)$ then $\Inn\pi_1(M)$ has elements of infinite order. 

Moreover, since $H^2(\pi_1(M),\Z)\cong H^2(M,\Z)$ and $H^1(\Zc\pi_1(M),\Z)\cong H^1(T^{\rank\Zc\pi_1(M)},\Z)$, the group $H^2(\Inn\pi_1(M),\Z)$ needs to be finitely generated. Using the results in \cite{adian2018central}, we can conclude that $\Inn\pi_1(M)$ cannot be isomorphic to a free Burnside group $B(a,b)$ with $b\geq 665$ odd. The results in \cite{chen2021actions} imply that $\Inn\pi_1(M)$ cannot be an infinite periodic 2-group of bounded exponent. Finally, we note that $\Inn\pi_1(M)$ is a finitely presented group and it is an open question if there exists finitely presented infinite periodic groups (this question is the Burnside problem for finitely presented groups, see \cite[pg. 3]{sapir2007some})


\selectlanguage{english}
\bibliographystyle{alpha}


\end{document}